\documentclass[10pt]{article}
\usepackage{amsfonts}
\usepackage{amsmath}
\usepackage{amscd,amssymb,amsthm}
\usepackage{graphics}
\usepackage{graphicx}

\newtheorem{theorem}{Theorem}[section]

\newtheorem{remark}[theorem]{Remark}

\begin{document}

\title{Analysis of a class of Lotka--Volterra systems }
\author{G. Moza\thanks{
Department of Mathematics, Politehnica University of Timisoara, Romania;
email: gheorghe.moza@upt.ro}, D. Constantinescu \thanks{
Department of Applied Mathematics, University of Craiova, Romania}, R.
Efrem, L. Bucur, R. Constantinescu }
\date{}
\maketitle

\begin{abstract}
A generalized two-dimensional cubic Lotka-Volterra model with infinitesimal
parameters is studied. Three different cases have been considered, one
non-degenerate and two degenerate. The local behavior of the model has been
studied in the three cases. Six bifurcation diagrams with thirty different
regions have been obtained in the non-degenerate case, respectively, sixteen
diagrams with forty regions in the two degenerate cases.
\end{abstract}

\section{Introduction}

Lotka Volterra systems are emblematic for the study of the interaction
between groups (populations) with different interests that coexist in the
same environment. Among the first models of this type are those presented in
[5] and [6]. While they are considered classical tools in the theory of
dynamical systems [7], [8], [9], [23], [24], [25], the interest for their
study has not decreased because they are successfully used for modeling
interactions in biology and ecology [10], [11], [12], [13], [14], but also
in economics [15], [16], chemistry [17] and engineering [18]. New variants
have been developed in recent years, for example, grey Lotka-Volterra or
fractional Lotka Volterra, and they have been proven effective in the study
of social or natural phenomena [19], [20], [21], [22] . All of these models
use intraspecific and interspecific interactions between species that are
expressed using first or second degree polynomials.

In this paper we propose and study a Lotka-Volterra two-dimensional system,
which generalizes a model published in \cite{tig1} and \cite{tig2}. More
exactly, we aim to study the system given by 
\begin{equation}
\left\{ 
\begin{tabular}{lll}
$\frac{dx}{d\tau }$ & $=$ & $2x\left[ \mu
_{1}+p_{11}x+p_{12}y+p_{13}xy+p_{14}x^{2}+p_{15}y^{2}\right] $ \\ 
&  &  \\ 
$\frac{dy}{d\tau }$ & $=$ & $2y\left[ \mu
_{2}+p_{21}x+p_{22}y+p_{23}xy+p_{24}x^{2}+p_{25}y^{2}\right] $%
\end{tabular}%
\right. ,  \label{g17}
\end{equation}%
where $p_{ij}=p_{ij}\left( \mu \right) $ are smooth functions of variable $%
\mu =\left( \mu _{1},\mu _{2}\right) \in 
\mathbb{R}
^{2}$ such that $p_{12}\left( 0\right) p_{21}\left( 0\right) \neq 0.$ The
present work is concerned with the study of the behavior of the system (\ref%
{g17}) when $\left\vert \mu _{1}\right\vert $ and $\left\vert \mu
_{2}\right\vert $ are infinitesimally small, that is, $\left\vert \mu
\right\vert =\sqrt{\mu _{1}^{2}+\mu _{2}^{2}}$ is sufficiently small, $%
0<\left\vert \mu \right\vert \ll 1.$ We approach in this article the case of
positive product of the nonzero parameters, that is, we assume $p_{12}\left(
0\right) p_{21}\left( 0\right) >0,$ and, more precisely, $p_{12}\left(
0\right) >0$ and $p_{21}\left( 0\right) >0.$ A similar model for the case $%
p_{12}\left( 0\right) p_{22}\left( 0\right) \neq 0$ is studied in \cite{tig1}
and \cite{tig2}.

\bigskip

In general, Lotka--Volterra models are defined for positive variables, $%
x\geq 0$ and $y\geq 0.$ In system (\ref{g17}), the axes $x=0$ and $y=0$ are
invariant with respect to the system's flow, thus, any orbit of this system
starting in the first quadrant remains in this region. So the first quadrant
is also an invariant region in the phase space of the system. As such, only
the behavior of the system (\ref{g17}) in this quadrant presents relevance
for practical applications. In this paper, we study the dynamics of the
system (\ref{g17}) in the first quadrant, where $x\geq 0$ and $y\geq 0.$

In order to decrease the number of parameters, we use the changes of
variables and time given by 
\begin{equation}
\xi _{1}=x\cdot p_{12}\left( \mu \right) \text{ and }\xi _{2}=y\cdot
p_{21}\left( \mu \right) ,t=2\tau .  \label{tr1}
\end{equation}%
By (\ref{tr1}), the system (\ref{g17}) is locally topologically equivalent
near the origin $O\left( 0,0\right) $ to 
\begin{equation}
\left\{ 
\begin{tabular}{lll}
$\dot{\xi}_{1}$ & $=$ & $\xi _{1}\left( \mu _{1}+\theta \left( \mu \right)
\xi _{1}+\gamma \left( \mu \right) \xi _{2}+M\left( \mu \right) \xi _{1}\xi
_{2}+N\left( \mu \right) \xi _{1}^{2}+L\left( \mu \right) \xi
_{2}^{2}\right) $ \\ 
&  &  \\ 
$\dot{\xi}_{2}$ & $=$ & $\xi _{2}\left( \mu _{2}+\frac{1}{\gamma \left( \mu
\right) }\xi _{1}+\delta \left( \mu \right) \xi _{2}+S\left( \mu \right) \xi
_{1}\xi _{2}+P\left( \mu \right) \xi _{2}^{2}+R\left( \mu \right) \xi
_{1}^{2}\right) $%
\end{tabular}%
\right. ,  \label{gl8}
\end{equation}%
where $\delta \left( \mu \right) =\frac{p_{22}\left( \mu \right) }{%
p_{21}\left( \mu \right) },$ $\theta \left( \mu \right) =\frac{p_{11}\left(
\mu \right) }{p_{12}\left( \mu \right) },$ $\gamma \left( \mu \right) =\frac{%
p_{12}\left( \mu \right) }{p_{21}\left( \mu \right) },$ $L\left( \mu \right)
=\frac{p_{15}\left( \mu \right) }{p_{21}^{2}\left( \mu \right) },$ $M\left(
\mu \right) =\frac{p_{13}\left( \mu \right) }{p_{12}\left( \mu \right)
p_{21}\left( \mu \right) },$ $N\left( \mu \right) =\frac{p_{14}\left( \mu
\right) }{p_{12}^{2}\left( \mu \right) },$ $P\left( \mu \right) =\frac{%
p_{25}\left( \mu \right) }{p_{21}^{2}\left( \mu \right) },$ $R\left( \mu
\right) =\frac{p_{24}\left( \mu \right) }{p_{12}^{2}\left( \mu \right) }$
and $S\left( \mu \right) =\frac{p_{23}\left( \mu \right) }{p_{12}\left( \mu
\right) p_{21}\left( \mu \right) }.$ The change $\left( x,y\right)
\longmapsto \left( \xi _{1},\xi _{2}\right) $ from (\ref{tr1}) is well
defined and nonsingular for all $\left\vert \mu \right\vert $ small enough,
because $p_{12}\left( 0,0\right) \neq 0$ and $p_{21}\left( 0,0\right) \neq
0; $ $\dot{\xi}_{1,2}=d\xi _{1,2}/dt.$

Since in what follows many of the coefficients are needed only at $\mu =0,$
we denote further by $\theta \left( 0\right) =\theta ,$ $\gamma \left(
0\right) =\gamma \neq 0,$ $\delta \left( 0\right) =\delta $ and so on; for
brevity, we denote $\left( 0,0\right) $ by $0.$

We study in this work the case of the positive product $p_{12}\left(
0\right) p_{21}\left( 0\right) >0,$ and consider $p_{12}\left( 0\right) >0$
and $p_{21}\left( 0\right) >0$ (mutualism case in Lotka-Volterra generalized
model). Thus, $\gamma >0$ while $\theta \in 
\mathbb{R}
$ and $\delta \in 
\mathbb{R}
.$

\begin{remark}
\label{rem1} Since $p_{12}\left( 0\right) >0$ and $p_{21}\left( 0\right) >0,$
the first quadrant of the system (\ref{g17}) is transformed by (\ref{tr1})
in the first quadrant $Q_{1}$ of the system (\ref{gl8}), which remains
invariant with respect to the flow of (\ref{gl8}). Thus, the new system (\ref%
{gl8}) will be studied in the first quadrant $Q_{1}$ defined by $\xi
_{1}\geq 0$ and $\xi _{2}\geq 0.$
\end{remark}

\begin{remark}
The equilibrium points are in the first quadrant $Q_{1}$ if their
coordinates are non-negative. We call them proper in this case, otherwise
virtual.
\end{remark}

\begin{remark}
\label{rem01} The system (\ref{g17}) with $p_{12}\left( 0\right) <0$ and $%
p_{21}\left( 0\right) <0,$ can also be reduced to a system of the form (\ref%
{gl8}), by changes $\xi _{1}=-x\cdot p_{12}\left( \mu \right) $ and $\xi
_{2}=-y\cdot p_{21}\left( \mu \right) ,$ respectively, $t=-2\tau .$ Indeed,
by these changes, (\ref{g17}) reduces to 
\begin{equation*}
\left\{ 
\begin{tabular}{lll}
$\dot{\xi}_{1}$ & $=$ & $\xi _{1}\left( -\mu _{1}+\theta \left( \mu \right)
\xi _{1}+\gamma \left( \mu \right) \xi _{2}-M\left( \mu \right) \xi _{1}\xi
_{2}-N\left( \mu \right) \xi _{1}^{2}-L\left( \mu \right) \xi
_{2}^{2}\right) $ \\ 
&  &  \\ 
$\dot{\xi}_{2}$ & $=$ & $\xi _{2}\left( -\mu _{2}+\frac{1}{\gamma \left( \mu
\right) }\xi _{1}+\delta \left( \mu \right) \xi _{2}-S\left( \mu \right) \xi
_{1}\xi _{2}-P\left( \mu \right) \xi _{2}^{2}-R\left( \mu \right) \xi
_{1}^{2}\right) $%
\end{tabular}%
\right. ,
\end{equation*}%
which is of the form (\ref{gl8}), if one denote the coefficients by $-\mu
_{1}=\nu _{1},$ $-\mu _{2}=\nu _{2},$ $-M\left( \mu \right) =M_{1}\left( \nu
\right) $ and so on. The new infinitesimal parameter is $\nu =\left( \nu
_{1},\nu _{2}\right) $ with $\left\vert \nu \right\vert =\left\vert \mu
\right\vert .$
\end{remark}

The paper is organized as follows. In Section 2 we study the behaviour of
the non-degenerate system, corresponding to $\theta \delta \neq 0.$ Section
3 is dedicated to the study of the degenerate system corresponding to $%
\delta =0$ and $\theta \neq 0,$ while the generate case $\theta =0$ and $%
\delta \neq 0$ is studied in Section 4. A summary of the obtained results
and some conclusions are presented in Section 5.

\section{The behavior of the system when $\protect\theta(0) \protect\delta%
(0) \neq 0$}

The first equilibrium of (\ref{gl8}) is $E_{0}\left( 0,0\right) .$ Two more
equilibria lying on the two axes, 
\begin{equation*}
E_{1}\left( -\frac{1}{\theta }\mu _{1}\left( 1+O\left( \left\vert \mu
\right\vert \right) \right) ,0\right) \text{ and }E_{2}\left( 0,-\frac{1}{%
\delta }\mu _{2}\left( 1+O\left( \left\vert \mu \right\vert \right) \right)
\right)
\end{equation*}%
bifurcate from $E_{0}$ as soon as $\mu _{1}\neq 0,$ respectively, $\mu
_{2}\neq 0,$ with $\left\vert \mu \right\vert $ small. $E_{1}$ is a proper
equilibrium if $\theta \mu _{1}<0,$ while $E_{2}$ is proper if $\delta \mu
_{2}<0;$ $\theta =\theta \left( 0\right) ,$ $\delta =\delta \left( 0\right) $
and $\gamma =\gamma \left( 0\right) .$

\begin{remark}
\label{rem3} The eigenvalues of the first three equilibria are $\mu _{1,2}$
of $E_{0},$ $-\mu _{1}+3\frac{N}{\theta ^{2}}\mu _{1}^{2}$ and $\mu _{2}-%
\frac{1}{\theta \gamma }\mu _{1}+\frac{R}{\theta ^{2}}\mu _{1}^{2}$ of $%
E_{1}\left( -\frac{1}{\theta }\mu _{1},0\right) ,$ respectively, $-\mu _{2}+3%
\frac{P}{\delta ^{2}}\mu _{2}^{2}$ and $\mu _{1}-\frac{\gamma }{\delta }\mu
_{2}+\frac{L}{\delta ^{2}}\mu _{2}^{2}$ of $E_{2}\left( 0,-\frac{1}{\delta }%
\mu _{2}\right) ,$ in their lowest terms.
\end{remark}

Another equilibrium $E_{3}\left( \xi _{1},\xi _{2}\right) $ is born close to 
$O$ for $\left\vert \mu \right\vert $ small, where 
\begin{equation*}
\xi _{1}=\left( \frac{-\delta \mu _{1}+\gamma \mu _{2}}{\theta \delta -1}%
\right) \left( 1+O\left( \left\vert \mu \right\vert \right) \right) \text{
and }\xi _{2}=\left( \frac{\mu _{1}-\theta \gamma \mu _{2}}{\gamma \left(
\theta \delta -1\right) }\right) \left( 1+O\left( \left\vert \mu \right\vert
\right) \right) ,
\end{equation*}%
provided that 
\begin{equation}
\theta \delta -1\neq 0.  \label{8a}
\end{equation}

The existence of $E_{3}$ is ensured by Implicit Functions Theorem applied to
the system 
\begin{equation}
\left\{ 
\begin{array}{c}
\mu _{1}+\theta \left( \mu \right) \xi _{1}+\gamma \left( \mu \right) \xi
_{2}+M\left( \mu \right) \xi _{1}\xi _{2}+N\left( \mu \right) \xi
_{1}^{2}+L\left( \mu \right) \xi _{2}^{2}=0 \\ 
\mu _{2}+\frac{1}{\gamma \left( \mu \right) }\xi _{1}+\delta \left( \mu
\right) \xi _{2}+S\left( \mu \right) \xi _{1}\xi _{2}+P\left( \mu \right)
\xi _{2}^{2}+R\left( \mu \right) \xi _{1}^{2}=0%
\end{array}%
\right. .  \label{ec2}
\end{equation}%
For $\left\vert \mu \right\vert $ sufficiently small, $E_{3}$ is proper when 
$(\mu _{1},\mu _{2})$ lies in the region 
\begin{equation*}
R=\left\{ \left( \mu _{1},\mu _{2}\right) \in 
\mathbb{R}
^{2}\mid \frac{-\delta \mu _{1}+\gamma \mu _{2}}{\delta \theta -1}>0,\frac{%
\mu _{1}-\theta \gamma \mu _{2}}{\gamma \left( \delta \theta -1\right) }%
>0\right\} .
\end{equation*}

\begin{theorem}
\label{th1} Assume $\left( \mu _{1},\mu _{2}\right) \in R.$ The following
assertions are true.

1) If $\theta \delta -1<0,$ then $E_{3}$ is a saddle.

2) If $\theta \delta -1>0,$ $E_{3}$ is an attractor (node or focus) when $%
\theta <0$ and $\delta <0,$ respectively, a repeller when $\theta >0$ and $%
\delta >0.$
\end{theorem}

\textit{Proof.} The characteristic polynomial at $E_{3}\left( \xi _{1},\xi
_{2}\right) $ is of the form $P\left( \lambda \right) =\lambda
^{2}-2p(\mu)\lambda +L(\mu)$ where

\begin{equation}
p(\mu )=\frac{1}{2}\left( \xi _{1}\theta (\mu )+\xi _{2}\delta (\mu )\right)
+\frac{1}{2}\left[ \xi _{1}\left( M(\mu )\xi _{2}+2N(\mu )\xi _{1}\right)
+\xi _{2}\left( 2P(\mu )\xi _{2}+S(\mu )\xi _{1}\right) \right]  \label{p1}
\end{equation}%
and

\begin{equation}
L(\mu )=\xi _{1}\xi _{2}\left( \theta (\mu )\delta (\mu )-1+c_{1}(\mu )\xi
_{1}+c_{2}(\mu )\xi _{2}+c_{3}(\mu )\xi _{1}^{2}+c_{4}(\mu )\xi _{1}\xi
_{2}+c_{5}(\mu )\xi _{2}^{2}\right) ,  \label{lgen}
\end{equation}%
$c_{1}(\mu )=2N(\mu )\delta (\mu )-\frac{1}{\gamma (\mu )}M(\mu )+S(\mu
)\theta (\mu )-2R(\mu )\gamma (\mu ),$

\noindent $c_{2}(\mu )=M(\mu )\delta (\mu )-S(\mu )\gamma (\mu )+2P(\mu
)\theta (\mu )-\frac{2L(\mu )}{\gamma (\mu )},$ $c_{3}(\mu )=-2(M(\mu )R(\mu
)-N(\mu )S(\mu )),$

\noindent $c_{4}(\mu )=-4(L(\mu)R(\mu )-N(\mu )P(\mu ))$ and $c_{5}(\mu
)=-2(L(\mu )S(\mu )-M(\mu )P(\mu)). $

Two bifurcation curves arise related to the existence of $E_{3},$ namely

\begin{equation*}
T_{1}=\left\{ \left( \mu _{1},\mu _{2}\right) \in 
\mathbb{R}
^{2}\left\vert \theta \gamma \mu _{2}=\mu _{1}+O\left( \mu _{1}^{2}\right)
\right. ,\theta \mu _{1}<0\right\}
\end{equation*}%
and 
\begin{equation*}
T_{2}=\left\{ \left( \mu _{1},\mu _{2}\right) \in 
\mathbb{R}
^{2}\left\vert \delta \mu _{1}=\gamma \mu _{2}+O\left( \mu _{2}^{2}\right)
\right. ,\delta \mu _{2}<0\right\} .
\end{equation*}

$E_{3}$ is born or vanishes when $\left( \mu _{1},\mu _{2}\right) $ crosses $%
T_{1}$ or $T_{2}.$ More exactly, $E_{1}\left( -\frac{\mu _{1}}{\theta }%
,0\right) $ collides to $E_{3}\left( -\frac{\mu _{1}}{\theta },0\right) $
when $\left( \mu _{1},\mu _{2}\right) \in T_{1},$ respectively, $E_{2}\left(
0,-\frac{1}{\delta }\mu _{2}\right) $ collides to $E_{3}\left( 0,-\frac{1}{%
\delta }\mu _{2}\right) $ on $T_{2};$ we call $E_{3}$ \textit{trivial} in
these cases, otherwise nontrivial.

It follows from (\ref{p1}) that 
\begin{equation}
p\left( \mu _{1},\mu _{2}\right) =\frac{1}{2\left( \theta \delta -1\right)
\gamma }\left[-\left(\theta\gamma-1\right)\delta\mu_{1}+\theta\gamma
\left(\gamma-\delta \right)\mu_{2} + O\left( \left\vert \mu\right\vert
^{2}\right) \right].  \label{p2}
\end{equation}%
Using Implicit Functions Theorem, the equation $p\left( \mu _{1},\mu
_{2}\right) =0$ exists as a unique curve in the parametric plane $\mu
_{1}\mu _{2}$ if $\theta \gamma -1\neq 0$ and $\gamma -\delta \neq 0,$ being
given by%
\begin{equation}
H=\left\{ \left( \mu _{1},\mu _{2}\right) \in 
\mathbb{R}
^{2}\left\vert \mu _{2}=\frac{\left( \theta \gamma -1\right) \delta }{\theta
\gamma \left( \gamma -\delta \right) }\mu _{1}+O\left( \mu _{1}^{2}\right)
\right. \right\} .  \label{has}
\end{equation}

1) Denoting the eigenvalues at $E_{3}$ by $\lambda _{1}$ and $\lambda _{2},$
it follows from $\theta \delta -1<0$ and (\ref{lgen}) that $\lambda
_{1}\lambda _{2}<0,$ whenever $E_{3}$ is nontrivial, thus, $E_{3}$ is a
saddle.

2) Assume further $\theta \delta -1>0$ and $\delta \neq \gamma .$ For $%
\left\vert \mu \right\vert $ sufficiently small, a nontrivial $E_{3}$ exists
in $Q_{1}$ iff $\left( \mu _{1},\mu _{2}\right) \in R_{1},$ where 
\begin{equation}
R_{1}=\left\{ \left. \left( \mu _{1},\mu _{2}\right) \in 
\mathbb{R}
^{2}\right\vert \text{ }\gamma \mu _{2}-\delta \mu _{1}>0,\mu _{1}-\theta
\gamma \mu _{2}>0,\left\vert \mu \right\vert <\varepsilon \right\} ,
\label{r1b}
\end{equation}%
where $0<\varepsilon \ll 1.$ Then $\lambda _{1}\lambda _{2}>0$ for all $%
\left( \mu _{1},\mu _{2}\right) \in R_{1},$ thus, the type of $E_{3}$ is
given by the sign of $p\left( \mu _{1},\mu _{2}\right) .$ We can obtain the
sign of $p$ by drawing the curve $H$ in the parametric plane $\mu _{1}\mu
_{2}$ and then determine the sign of $p$ on the two sides of $H.$
Alternatively, we observe from (\ref{has}) and (\ref{p2}) that only linear
terms in $\mu _{1}$ and $\mu _{2}$ are nedeed to define $H$ for $\left\vert
\mu \right\vert $ sufficiently small, which, in turn, implies that only
linear terms in $\xi _{1}$ and $\xi _{2}$ are used from (\ref{p1}) to obtain
the sign of $p.$ From $\delta \theta -1>0,$ it follows $\theta >0$ and $%
\delta >0$ or $\theta <0$ and $\delta <0.$ Thus, whenever $\left( \mu
_{1},\mu _{2}\right) \in R_{1},$ we obtain $p=\frac{1}{2}\left( \xi
_{1}\theta +\xi _{2}\delta \right) >0$ if $\theta >0$ and $\delta >0,$
respectively, $p=\frac{1}{2}\left( \xi _{1}\theta +\xi _{2}\delta \right) <0$
if $\theta <0$ and $\delta <0.$ Therefore, $\lambda _{1}+\lambda _{2}=2p>0$
in the first case, respectively, $\lambda _{1}+\lambda _{2}<0$ in the second
case, which, yield the conclusion. If $\delta \theta -1>0$ and $\delta
=\gamma >0,$ we have also $\theta >0$ and that implies $p>0.$ $\square $

\begin{remark}
\label{rem4} When $\theta \delta -1\neq 0,$ it follows from Theorem \ref{th1}
that the system (\ref{gl8}) does not undergo a Hopf bifurcation at $E_{3}.$
Indeed, $\lambda _{1},\lambda _{2}\in 
\mathbb{R}
$ if $\theta \delta -1<0,$ respectively, $p\left( \mu _{1},\mu _{2}\right)
\neq 0$ for $\left( \mu _{1},\mu _{2}\right) \in R_{1}$ if $\theta \delta
-1>0.$
\end{remark}

The analysis and results obtained so far in this section needed only terms
up to order two in $\xi _{1,2}$ from the system (\ref{gl8}), in a similar
way with the result obtained in (Lemma 8.17, \cite{Ku98}, page 361).
Therefore, the following conclusion can be drawn.

\begin{theorem}
\label{th2} If $p_{12}\left( 0\right) >0,$ $p_{21}\left( 0\right) >0$ and $%
\theta \delta -1\neq 0,$ the system (\ref{gl8}) is locally topologically
equivalent near the origin $O$ for all $\left\vert \mu \right\vert $
sufficiently small to the system%
\begin{equation}
\left\{ 
\begin{tabular}{lll}
$\frac{d\xi _{1}}{dt}$ & $=$ & $\xi _{1}\left( \mu _{1}+\theta \xi
_{1}+\gamma \xi _{2}\right) $ \\ 
$\frac{d\xi _{2}}{dt}$ & $=$ & $\xi _{2}\left( \mu _{2}+\frac{1}{\gamma }\xi
_{1}+\delta \xi _{2}\right) $%
\end{tabular}%
\right. .  \label{gl9}
\end{equation}
\end{theorem}

\begin{figure}[htpb]
\begin{center}
\includegraphics[width=0.3\textwidth]{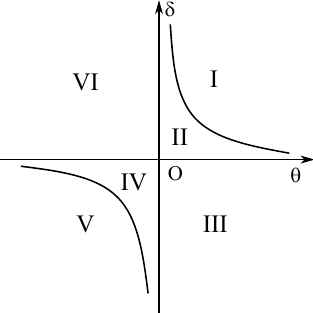}
\end{center}
\caption{When $\protect\theta\protect\delta\neq 0$ and $\protect\theta%
\protect\delta\neq 1,$ six cases in the $\protect\theta\protect\delta-$ plane
lead to six bifurcation diagrams.}
\label{fig1}
\end{figure}

\begin{figure}[htpb]
\begin{center}
\includegraphics[width=0.75\textwidth]{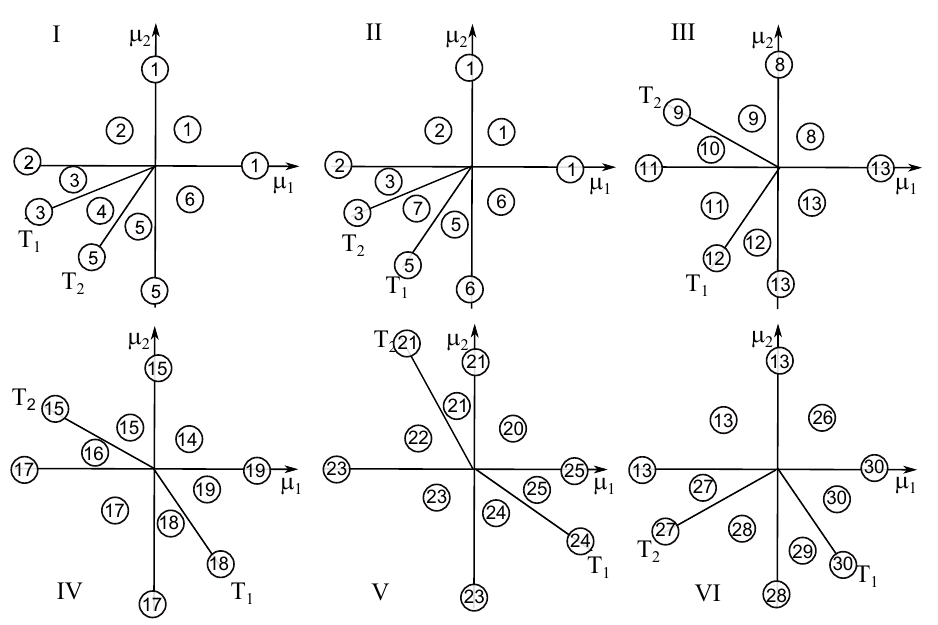}
\end{center}
\caption{Bifurcation diagrams for $\protect\theta\protect\delta\neq 0$
corresponding to the six regions I-VI.}
\label{fig2}
\end{figure}

In order to draw bifurcation diagrams, we notice that six cases arise in the 
$\theta \delta -$plane, which depend on the signs of $\theta ,$ $\delta $
and $\theta \delta -1,$ Figure \ref{fig1}. The cases give rise to six
bifurcation diagrams depicted in Figure \ref{fig2}, which contain 30
different regions in the parametric plane $\mu _{1}\mu _{2}.$ We describe in
Tables \ref{tabel1}--\ref{tabel2} the type of each equilibrium points from
the 30 regions, while the phase portraits corresponding to these regions are
depicted in Figures \ref{fig3}-\ref{fig4}.

\begin{table}[htpb]
\centering
\begin{tabular}{l|lllllllllllllll}
& $1$ & $2$ & $3$ & $4$ & $5$ & $6$ & $7$ & $8$ & $9$ & $10$ & $11$ & $12$ & 
$13$ & $14$ & $15$ \\ \hline
$E_{0}$ & $r$ & $s$ & $a$ & $a$ & $a$ & $s$ & $a$ & $r$ & $s$ & $s$ & $a$ & $%
a$ & $s$ & $r$ & $s$ \\ 
$E_{1}$ & $-$ & $r$ & $r$ & $s$ & $s$ & $-$ & $r$ & $-$ & $r$ & $r$ & $r$ & $%
s$ & $-$ & $s$ & $-$ \\ 
$E_{2}$ & $-$ & $-$ & $s$ & $s$ & $r$ & $r$ & $r$ & $s$ & $s$ & $a$ & $-$ & $%
-$ & $-$ & $s$ & $s$ \\ 
$E_{3}$ & $-$ & $-$ & $-$ & $r$ & $-$ & $-$ & $s$ & $-$ & $-$ & $s$ & $s$ & $%
-$ & $-$ & $-$ & $-$%
\end{tabular}%
\caption{\textit{The types of the equilibrium points of system \eqref{gl8}
for $\protect\theta\protect\delta\neq 0$ on different regions of the
bifurcation diagrams; the abbreviations s, a, r stand for saddle, attractor,
repeller, respectively.}}
\label{tabel1}
\end{table}

\begin{table}[htpb]
\centering
\begin{tabular}{l|lllllllllllllll}
& $16$ & $17$ & $18$ & $19$ & $20$ & $21$ & $22$ & $23$ & $24$ & $25$ & $26$
& $27$ & $28$ & $29$ & $30$ \\ \hline
$E_{0}$ & $s$ & $a$ & $s$ & $s$ & $r$ & $s$ & $s$ & $a$ & $s$ & $s$ & $r$ & $%
a$ & $a$ & $s$ & $s$ \\ 
$E_{1}$ & $-$ & $-$ & $a$ & $s$ & $s$ & $-$ & $-$ & $-$ & $a$ & $s$ & $s$ & $%
-$ & $-$ & $a$ & $s$ \\ 
$E_{2}$ & $a$ & $-$ & $-$ & $-$ & $s$ & $s$ & $a$ & $-$ & $-$ & $-$ & $-$ & $%
s$ & $r$ & $r$ & $r$ \\ 
$E_{3}$ & $s$ & $s$ & $s$ & $-$ & $a$ & $a$ & $-$ & $-$ & $-$ & $a$ & $-$ & $%
-$ & $s$ & $s$ & $-$%
\end{tabular}%
\caption{\textit{Continuation of Table \protect\ref{tabel1}.}}
\label{tabel2}
\end{table}

\begin{figure}[htpb]
\begin{center}
\includegraphics[width=1\textwidth]{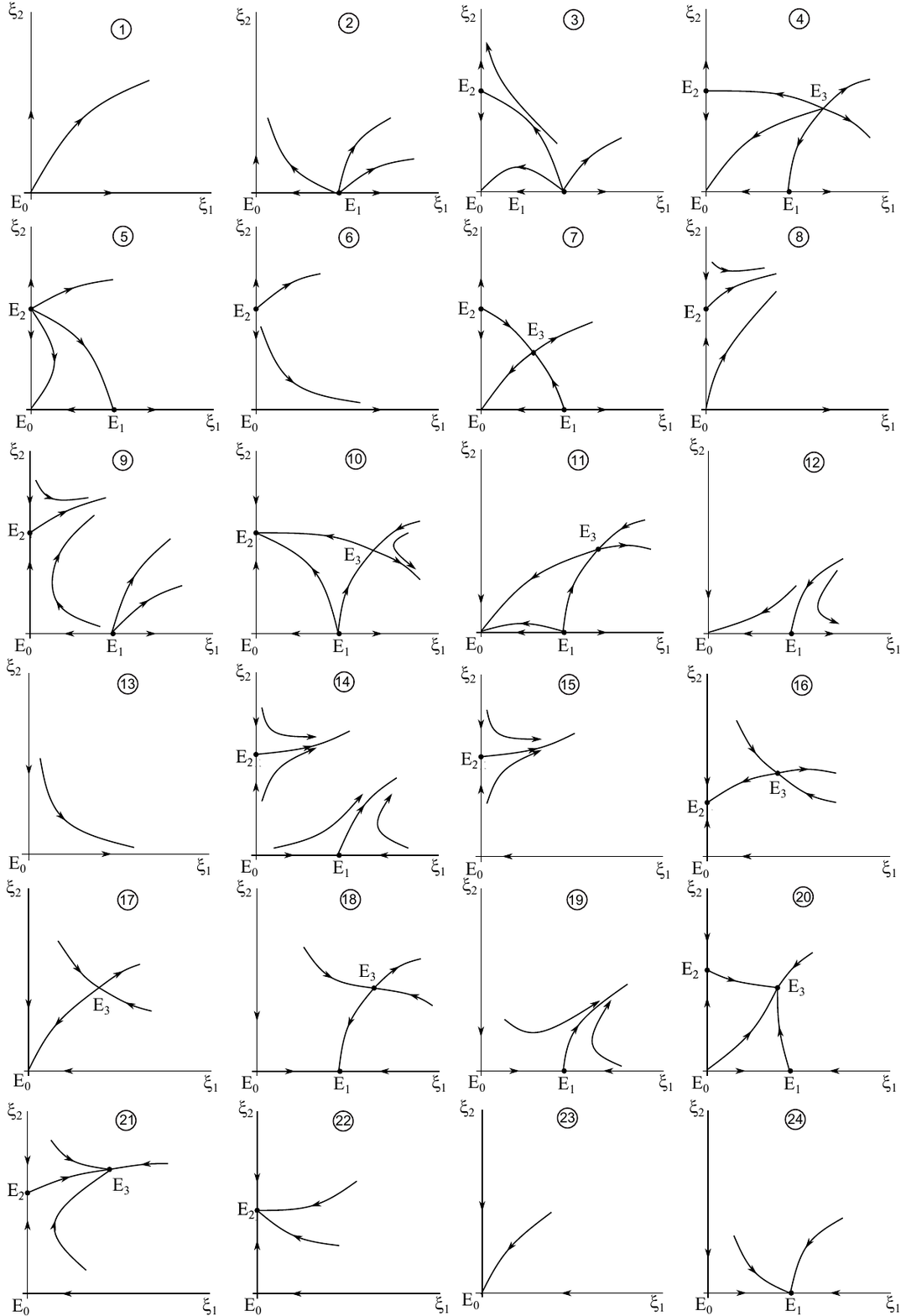}
\end{center}
\caption{Phase portraits corresponding to the bifurcation diagrams I-VI,
when $\protect\theta\protect\delta\neq 0.$}
\label{fig3}
\end{figure}

\begin{figure}[htpb]
\begin{center}
\includegraphics[width=1\textwidth]{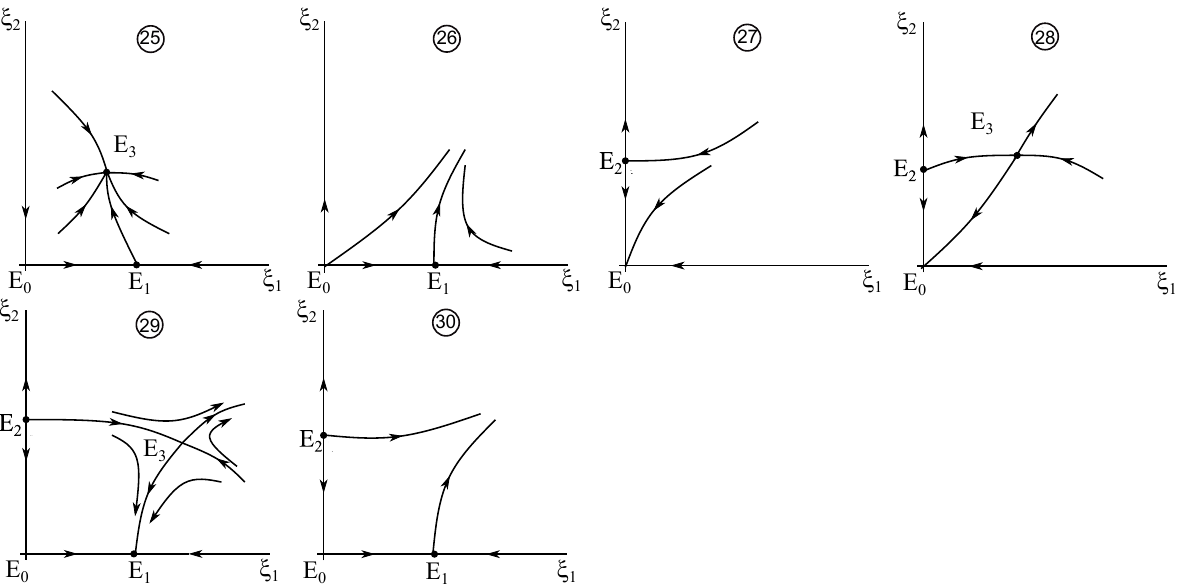}
\end{center}
\caption{Continuation of Figure \protect\ref{fig3}.}
\label{fig4}
\end{figure}

\section{The degeneracy of the system at $\protect\theta(0) \neq 0$ and $%
\protect\delta(0) =0$}

In the form (\ref{gl8}), the coefficients $\theta \left( \mu \right) ,$ $%
\delta \left( \mu \right) ,$ $\gamma \left( \mu \right) $ and the others are
smooth functions depending on the parameter $\mu =\left( \mu _{1},\mu
_{2}\right) .$ Thus, we can write 
\begin{equation*}
\delta \left( \mu \right) =\delta \left( 0\right) +\frac{\partial \delta
\left( 0\right) }{\partial \mu _{1}}\mu _{1}+\frac{\partial \delta \left(
0\right) }{\partial \mu _{2}}\mu _{2}+O\left( \left\vert \mu \right\vert
^{2}\right) ,
\end{equation*}%
and similarly for the other coefficients. Transforming (\ref{g17}) in (\ref%
{gl8}), a single constraint arises, namely $\gamma \left( 0\right) \neq 0,$
while $\theta \left( 0\right) $ and $\delta \left( 0\right) $ can be $0.$ In
the above section, we studied the case when neither of these constants is
zero, while here we consider $\theta \left( 0\right) =\theta \neq 0$ and $%
\delta \left( 0\right) =0.$ We keep $\gamma \left( 0\right) =\gamma >0$ and
assume $P\left( 0\right) \overset{not}{=}P\neq 0.$ Denote further by $\delta
_{i}=\frac{\partial \delta \left( 0\right) }{\partial \mu _{i}},$ $i=1,2,$
and assume $\delta _{1,2}\neq 0.$ \newline
Apart from $E_{0}\left( 0,0\right) $ and $E_{1}\left( -\frac{1}{\theta }\mu
_{1}\left( 1+O\left( \left\vert \mu \right\vert \right) \right) ,0\right) ,$
two more equilibrium points may exist, namely $E_{21}\left( 0,\xi
_{21}\right) $ and $E_{22}\left( 0,\xi _{22}\right) $ lying on the $\xi
_{2}- $axis, where 
\begin{equation*}
\xi _{21}=\frac{1}{2P\left( \mu \right) }\left( -\delta \left( \mu \right) +%
\sqrt{\Delta \left( \mu \right) }\right) \text{ and }\xi _{22}=\frac{1}{%
2P\left( \mu \right) }\left( -\delta \left( \mu \right) -\sqrt{\Delta \left(
\mu \right) }\right) ,
\end{equation*}%
whenever $\Delta \left( \mu \right) =\delta ^{2}\left( \mu \right) -4\mu
_{2}P\left( \mu \right) \geq 0.$ Notice that $\xi _{21}$ and $\xi _{22}$ are
the roots of $\mu _{2}+\delta \left( \mu \right) \xi _{2}+P\left( \mu
\right) \xi _{2}^{2}=0.$ In the lowest terms, $\Delta (\mu )$ reads $\Delta
(\mu )=\delta _{1}^{2}\mu _{1}^{2}\left( 1+O\left( \mu _{1}\right) \right)
-4P\mu _{2}\left( 1+O\left( \left\vert \mu \right\vert \right) \right) .$

For $\left\vert \mu \right\vert $ sufficiently small, denote by 
\begin{equation*}
D=\left\{ \left( \mu _{1},\text{ }\mu _{2}\right) \in 
\mathbb{R}
^{2}\left\vert \mu _{2}=\frac{\delta _{1}^{2}}{4P}\mu _{1}^{2}\left(
1+O\left( \mu _{1}\right) \right) \right. \right\} ,
\end{equation*}%
the bifurcation curve $\Delta \left( \mu \right) =0.$ The existence of $D$
is guaranteed by the Implicit Functions Theorem applied to $\Delta \left(
\mu _{1},\mu _{2}\right) =0,$ since $P=P(0)\neq 0.$ This curve is situated
on the first and the second quadrant if $P>0,$ respectively, third and
fourth quadrant if $P<0.$

\begin{remark}
\label{rem5} The eigenvalues of $E_{0}$ are $\mu _{1}$ and $\mu _{2},$ while
of $E_{1}$ they are $-\mu _{1}\left( 1+O\left( \left\vert \mu \right\vert
\right) \right) $ and $-\frac{1}{\theta \gamma }\left( \mu _{1}-\theta
\gamma \mu _{2}\right) +O(\left\vert \mu \right\vert ^{2}).$
\end{remark}

\begin{remark}
\label{r1} The eigenvalues of the equilibrium point $E_{2}(0,\xi _{2}),$
where $\xi _{2}$ satisfies $\mu _{2}+\delta (\mu )\xi _{2}+P(\mu )\xi
_{2}^{2}=0,$ are $\lambda _{1}^{E_{2}}=\xi _{2}\left( 2P(\mu )\xi
_{2}+\delta (\mu )\right) $ and $\lambda _{2}^{E_{2}}=L(\mu )\xi
_{2}^{2}+\gamma (\mu )\xi _{2}+\mu _{1}.$ Moreover $\lambda
_{1}^{E_{21}}=\xi _{21}\sqrt{\Delta (\mu )}$ and $\lambda _{1}^{E_{22}}=-\xi
_{22}\sqrt{\Delta (\mu )},$ thus, $\lambda _{1}^{E_{21}}>\lambda
_{1}^{E_{22}}.$
\end{remark}

The type of bifurcation by which $E_{11}$ and $E_{12}$ come into existence
or vanish, is described in the next theorem.

\begin{theorem}
\label{teo1} If $\theta \delta _{1}\delta _{2}P\left( 2P-\delta _{1}\gamma
\right) \neq 0$ and $\gamma >0,$ then $D\cap \left\{ (\mu _{1},\mu _{2})\mid
\mu _{1}<0\right\} $ and $D\cap \left\{ (\mu _{1},\mu _{2})\mid \mu
_{1}>0\right\} $ are two saddle-node bifurcation curves.
\end{theorem}

\textit{Proof.} Assume $\mu _{1}<0$ is fixed while $\mu _{2}$ varies, thus, $%
\mu _{2}$ is the bifurcation parameter. Write the system (\ref{gl8}) in the
form $\frac{d\xi }{dt}=f\left( \xi ,\mu \right) ,$ with $\xi =\left( \xi
_{1},\xi _{2}\right) ,$ $f=\left( f_{1},f_{2}\right) $ and $\mu =\left( \mu
_{1},\mu _{2}\right) .$ The proof is based on the Sotomayor's theorem \cite%
{Perko}. It is clear that $f\left( \xi _{0},\mu _{0}\right) =\left(
0,0\right) ,$ where $\xi _{0}=\left( 0,\xi _{21}\right) =\left( 0,\xi
_{22}\right) $ and $\mu _{0}=\left( \mu _{1},\mu _{2}\right) \in D\cap \{\mu
_{1}<0\}.$

On $D,$ the coinciding points $E_{21}\left( 0,\xi _{21}\right) $ and $%
E_{22}\left( 0,\xi _{22}\right) $ satisfy 
\begin{equation}
2P\left( \mu _{0}\right) \xi _{2}+\delta \left( \mu _{0}\right) =0,
\label{pe1}
\end{equation}%
where $\xi _{21}=\xi _{22}=\xi _{2}.$ Thus, the Jacobian matrix at $\left(
0,\xi _{2}\right) $ on $D$ becomes 
\begin{equation*}
A=Df\left( \xi _{0},\mu _{0}\right) =\left( 
\begin{array}{cc}
\mu _{1}+\gamma \left( \mu _{0}\right) \xi _{2}+L\left( \mu _{0}\right) \xi
_{2}^{2} & 0 \\ 
\frac{1}{\gamma \left( \mu _{0}\right) }\left( S\left( \mu _{0}\right)
\gamma \left( \mu _{0}\right) \xi _{2}+1\right) \xi _{2} & 0%
\end{array}%
\right) ,
\end{equation*}%
which has the eigenvector $v=\allowbreak \left( 
\begin{array}{cc}
0 & 1%
\end{array}%
\right) ^{T}$ corresponding to the eigenvalue $\lambda =0;$ $u^{T}$ denotes
the transpose of the vector $u.$ Similarly, $A^{T}$ has the eigenvalue $%
\lambda =0$ with the corresponding eigenvector $w.$ By (\ref{pe1}), $\xi
_{2}=\allowbreak -\frac{\delta _{1}}{2P}\mu _{1}\left( 1+O\left( \mu
_{1}\right) \right) $ in the expression of $Df\left( \xi _{0},\mu
_{0}\right) ,$ respectively, $w=\left( 
\begin{array}{cc}
\frac{\delta _{1}}{\gamma \left( 2P-\gamma \delta _{1}\right) }\left(
1+O\left( \mu _{1}\right) \right) & 1%
\end{array}%
\right) ^{T}.$

Then $f_{\mu _{2}}=\left( 
\begin{array}{c}
\frac{\partial f_{1}}{\partial \mu _{2}}\allowbreak \\ 
\frac{\partial f_{2}}{\partial \mu _{2}}\allowbreak%
\end{array}%
\right) =\left( 
\begin{array}{c}
0 \\ 
\xi _{2}\left( 1+\frac{\partial \delta }{\partial \mu _{2}}\xi _{2}+\frac{%
\partial P}{\partial \mu _{2}}\xi _{2}^{2}\right)%
\end{array}%
\right) $ at $\left( 0,\xi _{2}\right) ,$ thus, 
\begin{equation*}
w^{T}f_{\mu _{2}}\left( \xi _{0},\mu _{0}\right) =-\frac{\delta _{1}}{2P}\mu
_{1}\left( 1+O\left( \mu _{1}\right) \right) \neq 0.
\end{equation*}%
Further, determine $D^{2}f\left( \xi ,\mu \right) \left( v,v\right) =\left( 
\begin{array}{c}
d^{2}f_{1}\left( \xi ,\mu \right) \left( v,v\right) \\ 
d^{2}f_{2}\left( \xi ,\mu \right) \left( v,v\right)%
\end{array}%
\right) ,$ where 
\begin{equation*}
d^{2}f_{i}\left( \xi ,\mu \right) \left( v,v\right) =\frac{\partial ^{2}f_{i}%
}{\partial \xi _{1}^{2}}\left( \xi ,\mu \right) v_{1}^{2}+2\frac{\partial
^{2}f_{i}}{\partial \xi _{1}\partial \xi _{2}}\left( \xi ,\mu \right)
v_{1}v_{2}+\frac{\partial ^{2}f_{i}}{\partial \xi _{2}^{2}}\left( \xi ,\mu
\right) v_{2}^{2},i=1,2
\end{equation*}
is the differential of second order of the function $f_{i}$ applied to the
vector $v=\left( v_{1},v_{2}\right) =\left( 0,1\right)$. We obtain $%
d^{2}f_{1}\left( \xi _{0},\mu _{0}\right) \left( v,v\right) =0$ and $%
d^{2}f_{2}\left( \xi _{0},\mu _{0}\right) \left( v,v\right) =\allowbreak
2P\xi _{2}$ by (\ref{pe1}), thus,

\begin{equation*}
w^{T}\left[ D^{2}f\left( \xi _{0},\mu _{0}\right) \left( v,v\right) \right]
=-\delta _{1}\mu _{1}\left( 1+O\left( \mu _{1}\right) \right) \neq 0,
\end{equation*}%
which confirms the proof. For $\mu _{1}>0$ the proof is similar. $\square $

\begin{theorem}
\label{th3} If $\theta \delta _{1}\delta _{2}P\left( \delta _{1}\gamma
-P\right) \neq 0$ and $\gamma >0,$ the system (\ref{gl8}) has an equilibrium
point of the form $E_{3}\left( \xi _{1},\xi _{2}\right) ,$ which is a
saddle, whenever it is proper and non-trivial.
\end{theorem}

\textit{Proof.} From the Implicit Functions Theorem, the system (\ref{ec2})
has a unique solution of the form

\begin{equation}
\xi _{1}=\left( -\gamma \mu _{2}+\frac{1}{\gamma }\left( \delta _{1}\gamma
-P\right) \mu _{1}^{2}\right) \left( 1+O\left( \left\vert \mu \right\vert
\right) \right) \text{ and }\xi _{2}=\left( -\frac{1}{\gamma }\mu
_{1}+\theta \mu _{2}\right) (1+O\left( \left\vert \mu \right\vert \right) ),
\label{xi12}
\end{equation}%
for all $\left\vert \mu \right\vert $ sufficiently small. Thus, $E_{3}\left(
\xi _{1},\xi _{2}\right) $ with $\xi _{1,2}$ given by (\ref{xi12}) is an
equilibrium point of (\ref{gl8}). The characteristic polynomial at $E_{3}$
is $P\left( \lambda \right) =\lambda ^{2}-2p(\mu )\lambda +L(\mu )$ where 
\begin{equation}
L(\mu )=-\xi _{1}\xi _{2}\left( 1+O\left( \left\vert \xi \right\vert \right)
\right) .  \label{ell}
\end{equation}%
Thus, $E_{3}$ is a saddle whenever on $\xi _{1,2}>0,$ since its eigenvalues
satisfy $\lambda _{1}\lambda _{2}=L<0.$ $\square $

\bigskip

$E_{3}\left( \xi _{1},\xi _{2}\right) $ is well-defined when it lies in the
first quadrant given by $\xi _{1}\geq 0$ and $\xi _{2}\geq 0.$ Thus, two
bifurcation curves arise related to the existence of $E_{3},$ namely,

\begin{equation*}
T_{1}=\left\{ \left( \mu _{1},\mu _{2}\right) \in 
\mathbb{R}
^{2}\left\vert \mu _{2}=\frac{1}{\theta \gamma }\mu _{1}\left( 1+O\left( \mu
_{1}\right) \right) \right. ,\theta \mu _{1}<0\right\} ,
\end{equation*}%
given by $\xi _{2}=0,$ respectively, 
\begin{equation}
T_{3}=\left\{ \left. \left( \mu _{1},\mu _{2}\right) \in 
\mathbb{R}
^{2}\right\vert \mu _{2}=\frac{\delta _{1}\gamma -P}{\gamma ^{2}}\mu
_{1}^{2}\left( 1+O\left( \mu _{1}\right) \right) ,\text{ }\mu _{1}<0\right\}
,  \label{t3}
\end{equation}%
given by $\xi _{1}=0.$ $E_{3}$ collides to $E_{1}\left( -\frac{1}{\theta }%
\mu _{1}\left( 1+O\left( \mu _{1}\right) \right) ,0\right) $ on $T_{1},$
respectively, $E_{21}\left( 0,\xi _{21}\right) $ or $E_{22}\left( 0,\xi
_{22}\right) $ on $T_{3}.$ We call $E_{3}$ \textit{trivial} in these cases,
otherwise non-trivial.

In its lowest terms, $E_{3}$ reads 
\begin{equation*}
E_{3}\left( -\gamma \mu _{2}+\frac{1}{\gamma }\left( \delta _{1}\gamma
-P\right) \mu _{1}^{2},-\frac{1}{\gamma }\mu _{1}+\theta \mu _{2}\right) .
\end{equation*}

We assume further $P>0.$ A similar study can be performed for $P<0.$

\begin{theorem}
\label{th4} Assume $\theta \delta _{1}\neq 0,$ $P>0,$ $\gamma >0$ and $%
\left( \mu _{1},\mu _{2}\right) \in T_{3}$ with $\left\vert \mu \right\vert $
sufficiently small. Then,

a) if $\gamma \delta _{1}-2P<0,$ $E_{3}$ collides to $E_{21}\left( 0,-\frac{1%
}{\gamma }\mu _{1}\right) ,$ $\left. \lambda _{1}^{E_{21}}\right\vert
_{T_{3}}>0$ and $\left. \lambda _{2}^{E_{21}}\right\vert _{T_{3}}=0.$
Moreover, $E_{22}\left( 0,\mu _{1}\frac{P-\gamma \delta _{1}}{P\gamma }%
\right) $ is an attractor whenever it exists.

b) if $\gamma \delta _{1}-2P>0,$ $E_{3}$ collides to $E_{22}\left( 0,-\frac{1%
}{\gamma }\mu _{1}\right) ,$ $\left. \lambda _{1}^{E_{22}}\right\vert
_{T_{3}}<0$ and $\left. \lambda _{2}^{E_{22}}\right\vert _{T_{3}}=0.$
Moreover, $E_{21}\left( 0,\mu _{1}\frac{P-\gamma \delta _{1}}{P\gamma }%
\right) $ is a repeller whenever it exists.
\end{theorem}

\textit{Proof.} The eigenvalues of $E_{21}$ are $\lambda _{1}^{E_{21}}=-2\mu
_{2}-\mu _{1}\delta _{1}\xi _{21}$ and $\lambda _{2}^{E_{21}}=\mu
_{1}+\gamma \left( \mu \right) \xi _{21}+L\left( \mu \right) \xi _{21}^{2},$
while of $E_{22}$ they are $\lambda _{1}^{E_{22}}=-2\mu _{2}-\mu _{1}\delta
_{1}\xi _{22}$ and $\lambda _{2}^{E_{22}}=\mu _{1}+\gamma \left( \mu \right)
\xi _{22}+L\left( \mu \right) \xi _{22}^{2}.$ These lead to

\begin{equation}
\lambda _{1}^{E_{21}}=\xi _{21}\sqrt{\Delta }>0\text{ and }\lambda
_{1}^{E_{22}}=-\xi _{22}\sqrt{\Delta }<0.  \label{l11}
\end{equation}

Since it is difficult to study the signs of $\lambda _{2}^{E_{21}}$ and $%
\lambda _{2}^{E_{22}}$ as separate terms, we link them together through the
product 
\begin{equation}
\lambda _{2}^{E_{21}}\lambda _{2}^{E_{22}}=\frac{1}{P}\left( \mu _{2}\gamma
^{2}-\left( \delta _{1}\gamma -P\right) \mu _{1}^{2}\right) \left( 1+O\left(
\left\vert \mu \right\vert \right) \right) .  \label{l22}
\end{equation}%
a) Assume $\gamma \delta _{1}-2P<0.$ Then $\left. \Delta \left( \mu \right)
\right\vert _{T_{3}}=\frac{\mu _{1}^{2}}{\gamma ^{2}}\left( \gamma \delta
_{1}-2P\right) ^{2}$ leads to $\xi _{21}=\allowbreak \allowbreak -\frac{1}{%
\gamma }\mu _{1}$ and $\xi _{22}=\mu _{1}\frac{P-\gamma \delta _{1}}{P\gamma 
},$ where $\mu _{1}<0.$ Thus, $E_{3}$ collides to $E_{21}\left( 0,-\frac{1}{%
\gamma }\mu _{1}\right) \ $for $\left\vert \mu \right\vert $ sufficiently
small and $\left( \mu _{1},\mu _{2}\right) \in T_{3}.$ Notice that $%
E_{3}\left( 0,\xi _{2}\right) $ satisfies 
\begin{equation}
\mu _{1}+\gamma \left( \mu \right) \xi _{2}+L\left( \mu \right) \xi
_{2}^{2}=0\text{ and }\mu _{2}+\delta \left( \mu \right) \xi _{2}+P\left(
\mu \right) \xi _{2}^{2}=0,  \label{rel1}
\end{equation}%
thus, the eigenvalues of the colliding points $E_{3}$ and $E_{21}$ on $T_{3}$
satisfy $\left. \lambda _{1}^{E_{21}}\right\vert _{T_{3}}=-\frac{\gamma
\delta _{1}-2P}{\gamma ^{2}}\mu _{1}^{2}>0$ and $\left. \lambda
_{2}^{E_{21}}\right\vert _{T_{3}}=0.$

In addition, $E_{22}\left( 0,\mu _{1}\frac{P-\gamma \delta _{1}}{P\gamma }%
\right) $ has its eigenvalues $\left. \lambda _{1}^{E_{22}}\right\vert
_{T_{3}}=\left. -\xi _{22}\sqrt{\Delta }\right\vert _{T_{3}}<0$ and 
\begin{equation*}
\lambda _{2}^{E_{22}}|_{T_{3}}=-\mu _{1}\frac{\gamma \delta _{1}-2P}{P}%
\left( 1+O\left( \mu _{1}\right) \right) <0,
\end{equation*}%
since $\mu _{1}<0$ on $T_{3}.$ Thus, $E_{22}$ is an attractor on $T_{3}.$

Since $\lambda _{2}^{E_{22}}|_{T_{3}}\neq 0$ on $T_{3},$ the curve $\lambda
_{2}^{E_{21}}\lambda _{2}^{E_{22}}=0$ given by (\ref{l22}) coincides to $%
T_{3}$ and $\lambda _{2}^{E_{21}}=0.$

b) Assume $\gamma \delta _{1}-2P>0.$ Then, $E_{3}$ collides to $E_{22}\left(
0,-\frac{1}{\gamma }\mu _{1}\right) \ $for $\left\vert \mu \right\vert $
sufficiently small and $\mu \in T_{3},$ with the eigenvalues $\left. \lambda
_{1}^{E_{22}}\right\vert _{T_{3}}=\left. -\xi _{22}\sqrt{\Delta }\right\vert
_{T_{3}}<0$ and $\left. \lambda _{2}^{E_{22}}\right\vert _{T_{3}}=0,$ by (%
\ref{rel1}). Also, the eigenvalues of $E_{21}\left( 0,\mu _{1}\frac{P-\gamma
\delta _{1}}{P\gamma }\right) $ are $\left. \lambda _{1}^{E_{21}}\right\vert
_{T_{3}}=\xi _{21}\sqrt{\Delta }>0$ and 
\begin{equation*}
\lambda _{2}^{E_{21}}|_{T_{3}}=-\mu _{1}\frac{\gamma \delta _{1}-2P}{P}%
\left( 1+O\left( \mu _{1}\right) \right) >0,
\end{equation*}%
since $\mu _{1}<0$ on $T_{3}.$ Thus, $E_{21}$ is a repeller on $T_{3}.$

Since $\lambda _{2}^{E_{21}}|_{T_{3}}\neq 0$ on $T_{3},$ the curve $\lambda
_{2}^{E_{21}}\lambda _{2}^{E_{22}}=0$ given by (\ref{l22}) coincides to $%
T_{3}$ and $\lambda _{2}^{E_{22}}=0.$ $\square $

\begin{remark}
The results obtained in Theorem \ref{th4} are important because they show
that the behaviors of $E_{21}$ and $E_{22}$ are determined by the
bifurcation curve $T_{3},$ which is related to another equilibrium point,
namely $E_{3}.$\newline
\end{remark}

\begin{theorem}
Assume $\delta _{1}\gamma _{2}\left( \delta _{1}\gamma -2P\right) \neq 0,$ $%
P>0$ and $\gamma >0,$ where $\gamma _{2}=\frac{\partial \gamma }{\partial
\mu _{2}}(0).$ Then $T_{3}$ is a transcritical bifurcation curve.
\end{theorem}

\textit{Proof.} Assume first $\gamma \delta _{1}-2P<0,$ thus, $E_{3}$
collides to $E_{21}\left( 0,-\frac{1}{\gamma }\mu _{1}\right) $ on $T_{3}.$
For $\gamma \delta _{1}-2P>0,$ the proof is similar. Let $\mu _{2}$ be the
bifurcation parameter while $\mu _{1}<0$ is assumed fixed. Denote by $\mu
_{0}=\left( \mu _{1},\mu _{2}\right) \in T_{3},$ that is, $\mu _{2}=\frac{%
\delta _{1}\gamma -P}{\gamma ^{2}}\mu _{1}^{2}.$

On $T_{3},$ $E_{3}\left( 0,\xi _{2}\right) $ satisfies (\ref{rel1}), thus, $%
\xi _{2}=\xi _{21}=-\frac{\mu _{1}}{\gamma \left( \mu _{0}\right) }.$ For $%
\left( \mu _{1},\mu _{2}\right) \in T_{3},$ denote by $\xi _{0}=\left( 0,\xi
_{2}\right) .$ The Jacobian matrix at $\left( \xi _{0},\mu _{0}\right) $
reads

\begin{equation*}
A=Df\left( \xi _{0},\mu _{0}\right) =\allowbreak \left( 
\begin{array}{cc}
0 & 0 \\ 
\frac{1}{\gamma \left( \mu _{0}\right) }\xi _{2}\left( S\left( \mu
_{0}\right) \gamma \left( \mu _{0}\right) \xi _{2}+1\right) & \delta \left(
\mu _{0}\right) \xi _{2}+2P\left( \mu _{0}\right) \xi _{2}^{2}%
\end{array}%
\allowbreak \right) ,
\end{equation*}%
where $f=\left( f_{1},f_{2}\right) .$ Then $A$ and $A^{T}$ have the
eigenvalue $0$ with the corresponding eigenvector $v=\left( 
\begin{array}{cc}
v_{1} & 1%
\end{array}%
\right) ^{T},$ $v_{1}=$ $-\frac{\gamma \left( \mu _{0}\right) \delta \left(
\mu _{0}\right) +2P\left( \mu _{0}\right) \gamma \left( \mu _{0}\right) \xi
_{2}}{S\left( \mu _{0}\right) \gamma \left( \mu _{0}\right) \xi _{2}+1},$
for $A,$ respectively, $w=\allowbreak \left( 
\begin{array}{cc}
1 & 0%
\end{array}%
\right) ^{T}$ for $A^{T}.$

Further, $f_{\mu _{2}}$ has the form

\begin{equation*}
f_{\mu _{2}}=\left( 
\begin{array}{c}
\frac{\partial f_{1}}{\partial \mu _{2}} \\ 
\frac{\partial f_{2}}{\partial \mu _{2}}%
\end{array}%
\right) =\left( 
\begin{array}{c}
\xi _{1}\left( \frac{\partial \theta \left( \mu \right) }{\partial \mu _{2}}%
\xi _{1}+\frac{\partial \gamma \left( \mu \right) }{\partial \mu _{2}}\xi
_{2}+\frac{\partial M\left( \mu \right) }{\partial \mu _{2}}\xi _{1}\xi _{2}+%
\frac{\partial N\left( \mu \right) }{\partial \mu _{2}}\xi _{1}^{2}+\frac{%
\partial L\left( \mu \right) }{\partial \mu _{2}}\xi _{2}^{2}\right) \\ 
\xi _{2}\left( 1+B\left( \mu \right) \xi _{1}+\frac{\partial \delta \left(
\mu \right) }{\partial \mu _{2}}\xi _{2}+\frac{\partial S\left( \mu \right) 
}{\partial \mu _{2}}\xi _{1}\xi _{2}+\frac{\partial P\left( \mu \right) }{%
\partial \mu _{2}}\xi _{2}^{2}+\frac{\partial R\left( \mu \right) }{\partial
\mu _{2}}\xi _{1}^{2}\right)%
\end{array}%
\right) ,
\end{equation*}%
where $B\left( \mu \right) =-\frac{1}{\gamma ^{2}\left( \mu \right) }\frac{%
\partial \gamma \left( \mu \right) }{\partial \mu _{2}}.$ Then $%
C_{1}=w^{T}\cdot f_{\mu _{2}}\left( \xi _{0},\mu _{0}\right) =0.$

The Jacobian $Df_{\mu _{2}}$ at the point $\left( \xi _{0},\mu _{0}\right) $
applied to the vector $v$ has the form

\begin{equation*}
Df_{\mu _{2}}\left( \xi _{0},\mu _{0}\right) \left( v\right) =\left( 
\begin{array}{cc}
\frac{\partial \gamma \left( \mu _{0}\right) }{\partial \mu _{2}}\xi _{2}+%
\frac{\partial L\left( \mu _{0}\right) }{\partial \mu _{2}}\xi _{2}^{2} & 0
\\ 
\frac{\partial ^{2}f_{2}}{\partial \xi _{1}\partial \mu _{2}}\left( \xi
_{0},\mu _{0}\right) & \frac{\partial ^{2}f_{2}}{\partial \xi _{2}\partial
\mu _{2}}\left( \xi _{0},\mu _{0}\right)%
\end{array}%
\allowbreak \right) \left( 
\begin{array}{c}
v_{1} \\ 
1%
\end{array}%
\right) =\left( 
\begin{array}{c}
v_{1}\left( \frac{\partial \gamma \left( \mu _{0}\right) }{\partial \mu _{2}}%
\xi _{2}+\frac{\partial L\left( \mu _{0}\right) }{\partial \mu _{2}}\xi
_{2}^{2}\right) \\ 
K%
\end{array}%
\right) ,
\end{equation*}%
where $K$ is an expression which is not needed in what follows. So the
second coefficient is 
\begin{equation*}
C_{2}=w^{T}\cdot \left[ Df_{\mu _{2}}\left( \xi _{0},\mu _{0}\right) \left(
v\right) \right] =v_{1}\left( \frac{\partial \gamma \left( \mu _{0}\right) }{%
\partial \mu _{2}}\xi _{2}+\frac{\partial L\left( \mu _{0}\right) }{\partial
\mu _{2}}\xi _{2}^{2}\right) =\mu _{1}^{2}\frac{\left( \gamma \delta
_{1}-2P\right) \gamma _{2}}{\gamma }\left( 1+O\left( \mu _{1}\right) \right)
\neq 0.
\end{equation*}

It remains to find $C_{3}=w^{T}\left[ D^{2}f\left( \xi _{0},\mu _{0}\right)
\left( v,v\right) \right] ,$ where $D^{2}f\left( \xi ,\mu \right) \left(
v,v\right) =\left( 
\begin{array}{c}
d^{2}f_{1}\left( \xi ,\mu \right) \left( v,v\right) \\ 
d^{2}f_{2}\left( \xi ,\mu \right) \left( v,v\right)%
\end{array}%
\right) .$ Since $w=\allowbreak \left( 
\begin{array}{cc}
1 & 0%
\end{array}%
\right) ^{T},$ only $d^{2}f_{1}$ is needed. We obtain

\begin{equation*}
C_{3}=2\gamma \mu _{1}\left( 2P-\gamma \delta _{1}\right) \left( 1+O\left(
\mu _{1}\right) \right) \neq 0.
\end{equation*}

\begin{remark}
One can show similarly that $T_{3},$ $X^{+}=\left\{ \left( \mu _{1},0\right)
\left\vert \mu _{1}>0\right. \right\} ,$ $X^{-}=\left\{ \left( \mu
_{1},0\right) \left\vert \mu _{1}<0\right. \right\} ,$ $Y^{+}=\left\{ \left(
0,\mu _{2}\right) \left\vert \mu _{2}>0\right. \right\} $ and $Y^{-}=\left\{
\left( 0,\mu _{2}\right) \left\vert \mu _{2}<0\right. \right\} $ are
transcritical bifurcation curves.
\end{remark}

\begin{remark}
The behavior of the system (\ref{gl8}) on the axes $X^{\pm }$ and $Y^{\pm },$
and on the transcritical curves $T_{1}$ and $T_{3},$  coincides with the
behavior of (\ref{gl8}) on the left or right regions delimited by these
curves, where one of the collinding points became virtual after collision.
See also \cite{tig1} and \cite{tig2}. On the saddle node curves, the
corresponding dynamics is presented in Fig.\ref{fig8}.
\end{remark}

Restricted to quadrant I, the phase portraits on the bifurcation curves
coincide to the phase portraits corresponding to the regions where one of
the collinding points became virtual after collision.

\bigskip 

Define the following regions

\begin{equation*}
R_{00}=\left\{ \left( \mu _{1},\mu _{2}\right) \in 
\mathbb{R}
^{2}\left\vert \Delta \left( \mu \right) >0\right. ,\mu _{2}>0,\delta
_{1}\mu _{1}>0\right\} \cup \left\{ \left( \mu _{1},\mu _{2}\right) \in 
\mathbb{R}
^{2}\left\vert \Delta \left( \mu \right) <0\right. \right\} ,
\end{equation*}

\begin{equation*}
R_{10}=\left\{ \left( \mu _{1},\mu _{2}\right) \in 
\mathbb{R}
^{2}\left\vert \mu _{2}<0\right. \right\} \text{ and }R_{20}=\left\{ \left(
\mu _{1},\mu _{2}\right) \in 
\mathbb{R}
^{2}\left\vert \Delta \left( \mu \right) >0\right. ,\mu _{2}>0,\delta
_{1}\mu _{1}<0\right\} .
\end{equation*}

Both points $E_{21}$ and $E_{22}$ are proper in the region $R_{20},$ while $%
R_{10}$ contains only $E_{21}$ proper, since $\xi _{21}>\xi _{22}$ and $P>0.$
On $R_{00},$ none of the two equilibria survive (they are virtual), because
either $\xi _{21}<0$ and $\xi _{22}<0$ or $\xi _{21}$ and $\xi _{22}$ are
not real numbers. Let us show that $\xi _{21}>0$ and $\xi _{22}>0$ on $%
R_{20}.$ To this end, we show $\xi _{21}\xi _{22}>0$ and $\xi _{21}+\xi
_{22}>0.$ Since, in their lowest terms, $\xi _{21}\xi _{22}=\frac{\mu _{2}}{P%
}>0$ on $R_{20},$ it remains to prove $\xi _{21}+\xi _{22}=-\frac{1}{P}%
\left( \delta _{1}\mu _{1}+\delta _{2}\mu _{2}\right) >0.$ It is clear that $%
\xi _{21}+\xi _{22}>0$ on $R_{20}$ if $\delta _{2}<0.$ Assume $\delta
_{2}>0. $ From $\Delta \left( \mu \right) >0\allowbreak $ we get $\mu _{2}<%
\frac{\mu _{1}^{2}\delta _{1}^{2}}{4P},$ which, in turns, implies 
\begin{equation*}
\xi _{21}+\xi _{22}>\frac{-\mu _{1}\delta _{1}}{P}\left( 1+O\left( \mu
_{1}\right) \right) >0.
\end{equation*}

\begin{remark}
If $\left( \delta _{1}\gamma -2P\right) P\neq 0,$ then the parabola $T_{3}$
is situated under the parabola $D=\left\{ \left( \mu _{1},\mu _{2}\right)
\in 
\mathbb{R}
^{2}\mid \mu _{2}=\frac{\delta _{1}^{2}\mu _{1}^{2}}{4P}(1+O\left( \mu
_{1}\right) \right\} .$
\end{remark}

Whenever $\delta _{1}\gamma -P>0,$ denote by 
\begin{equation*}
R_{20}^{-}=R_{20}\cap \left\{ \mu _{2}\gamma ^{2}<\left( \delta _{1}\gamma
-P\right) \mu _{1}^{2}\right\} \text{ and }R_{20}^{+}=R_{20}\cap \left\{ \mu
_{2}\gamma ^{2}>\left( \delta _{1}\gamma -P\right) \mu _{1}^{2}\right\} ,
\end{equation*}%
the regions from $R_{20}$ to the left, respectively, the right of $T_{3}.$
Notice that 
\begin{equation*}
R_{20}=R_{20}^{-}\cup T_{3}\cup R_{20}^{+}.
\end{equation*}

\begin{theorem}
\label{th5} Assume $\theta \delta _{1}\neq 0,\gamma >0,$ $P>0$ and $\gamma
\delta _{1}-P>0.$ Then,

a) if $\gamma \delta _{1}-2P<0,$ then $E_{21}$ is a repeller on $R_{10}\cup
X_{-}\cup R_{20}^{-}$ and a saddle on $R_{20}^{+},$ while $E_{22}$ is an
attractor on $R_{20}.$

b) if $\gamma \delta _{1}-2P>0,$ then $E_{21}$ is a repeller on $R_{10}\cup
X_{-}\cup R_{20},$ while $E_{22}$ is a saddle on $R_{20}^{+}$ and an
attractor on $R_{20}^{-}.$
\end{theorem}

\textit{Proof.} a) Notice first that $\delta _{1}>0$ and $0<P<\gamma \delta
_{1}<2P,$ thus, $R_{20}\subset \left\{ \left( \mu _{1},\mu _{2}\right) \in 
\mathbb{R}
^{2}\left\vert \mu _{1}<0\right. ,\mu _{2}>0\right\} .$ Also $T_{3}$ and $%
D\cap \{\mu _{1}<0\}$ are included in the second quadrant. From Theorem \ref%
{th4}, $E_{3}$ collides to $E_{21}\left( 0,-\frac{1}{\gamma }\mu _{1}\right) 
$ on $T_{3}$ and $\left. \lambda _{2}^{E_{21}}\right\vert _{T_{3}}=0.$ Thus, 
$\lambda _{2}^{E_{21}}\left( \mu _{1},\mu _{2}\right) $ keeps constant sign
outside $T_{3}$ and changes its sign when $\left( \mu _{1},\mu _{2}\right) $
crosses $T_{3}.$ Since 
\begin{equation*}
\lambda _{2}^{E_{21}}\left( 0,\mu _{2}\right) =\frac{\gamma }{P}\sqrt{-\mu
_{2}P}\left( 1+O\left( \mu _{2}\right) \right) >0
\end{equation*}
if $\mu _{2}<0,$ it follows that $\lambda _{2}^{E_{21}}>0$ on $R_{20}^{-}$
and $\lambda _{2}^{E_{21}}<0$ on $R_{20}^{+},$ while $\lambda
_{1}^{E_{21}}=\xi _{21}\sqrt{\Delta }>0$ on $R_{20}.$ Thus, $E_{21}$ is a
saddle on $R_{20}^{+}$ and a repeller on $R_{20}^{-}.$ It remains a repeller
on $R_{10}\cup X_{-}$ based on the same reasons.

On the other hand, since $\left. \lambda _{2}^{E_{22}}\right\vert
_{T_{3}}\neq 0,$ $\lambda _{2}^{E_{22}}$ keeps constant sign $R_{20},$ which
is negative by Theorem \ref{th4}. From $\lambda _{1}^{E_{22}}=-\xi _{22}%
\sqrt{\Delta }<0,$ it follows that $E_{22}$ is an attractor on $R_{20}.$

b) We have in this case $\delta _{1}>0$ and $0<2P<\gamma \delta _{1},$ which
yield $R_{20}\subset \left\{ \left( \mu _{1},\mu _{2}\right) \in 
\mathbb{R}
^{2}\left\vert \mu _{1}<0\right. ,\mu _{2}>0\right\} .$ From Theorem \ref%
{th4}, $E_{3}$ collides to $E_{22}\left( 0,-\frac{1}{\gamma }\mu _{1}\right) 
$ on $T_{3}$ and $\left. \lambda _{2}^{E_{22}}\right\vert _{T_{3}}=0,$ thus, 
$\lambda _{2}^{E_{22}}\left( \mu _{1},\mu _{2}\right) $ changes its sign
when $\left( \mu _{1},\mu _{2}\right) $ crosses $T_{3}.$More exactly, $%
\lambda _{2}^{E_{22}}>0$ on $R_{20}^{+}$ because 
\begin{equation*}
\left. \lambda _{2}^{E_{22}}\right\vert _{D}=-\frac{1}{2P}\mu _{1}\left(
\gamma \delta _{1}-2P\right) \left( 1+O\left( \mu _{1}\right) \right) >0,
\end{equation*}
which, in turn, yields $\lambda _{2}^{E_{22}}<0$ on $R_{20}^{-}.$ Thus, $%
E_{22}$ is a saddle on $R_{20}^{+}$ and an attractor on $R_{20}^{-},$
because $\lambda _{1}^{E_{22}}<0$ on $R_{20}.$

Further, $\lambda _{2}^{E_{21}}\neq 0$ on $R_{10}\cup X_{-}\cup R_{20}$ and $%
\left. \lambda _{2}^{E_{21}}\right\vert _{T_{3}}=-\mu _{1}\frac{\gamma
\delta _{1}-2P}{P}\left(1+O\left(\mu_{1}\right)\right)>0,$ yield $\lambda
_{2}^{E_{21}}>0$ on $R_{10}\cup X_{-}\cup R_{20}.$ Thus, $E_{21}$ is a
repeller on $R_{10}\cup X_{-}\cup R_{20}.$ $\square $\newline

\bigskip

Whenever $\delta _{1}\gamma -P<0,$ denote by $R_{10}^{+}=R_{10}\cap \left\{
\left( \mu _{1},\mu _{2}\right) \in 
\mathbb{R}
^{2}\left\vert \mu _{2}\gamma ^{2}>\left( \delta _{1}\gamma -P\right) \mu
_{1}^{2}\right. ,\mu _{1}<0\right\} $ and $R_{10}^{-}\subset R_{10}$ such
that $R_{10}=R_{10}^{+}\cup T_{3}\cup R_{10}^{-}.$ Denote by 
\begin{equation*}
T_{3}^{+}=\left\{ \left( \mu _{1},\mu _{2}\right) \in 
\mathbb{R}
^{2}\left\vert \gamma ^{2}\mu _{2}=\left( \delta _{1}\gamma -P\right) \mu
_{1}^{2}\left( 1+O\left( \mu _{1}\right) \right) \right. ,\text{ }\mu
_{1}>0\right\} .
\end{equation*}

\begin{theorem}
\label{th6} Assume $\theta \delta _{1}\neq 0,\gamma >0,$ $P>0$ and $\gamma
\delta _{1}-P<0.$ Then,

a) if $\delta _{1}>0,$ $E_{21}$ is a saddle and $E_{22}$ an attractor on $%
R_{20}.$ Moreover, $E_{21}$ is a saddle on $X_{-}\cup R_{10}^{+}$ and a
repeller on $R_{10}^{-}.$

b) if $\delta _{1}<0,$ $E_{21}$ is a repeller and $E_{22}$ a saddle on $%
R_{20}.$ Moreover, $E_{21}$ is a repeller on $X_{+}\cup R_{10}^{-}$ and a
saddle on $R_{10}^{+}.$
\end{theorem}

\textit{Proof.} Notice that only $\gamma \delta _{1}-2P<0$ is possible in
this case, while $\gamma \delta _{1}-2P>0$ leads to a contradiction.

a) Assume first $\delta _{1}>0,$ thus, $R_{20}\subset \left\{ \left( \mu
_{1},\mu _{2}\right) \in 
\mathbb{R}
^{2}\left\vert \mu _{1}<0\right. ,\mu _{2}>0\right\} $ and $T_{3}\subset
R_{10}.$ If $\mu _{1}<0,$ then $\lambda _{2}^{E_{21}}\lambda _{2}^{E_{22}}=0$
only on $T_{3}\nsubseteq R_{20}\cup X_{-}\cup R_{10}^{+}.$ Thus, $\lambda
_{2}^{E_{21}}$ and $\lambda _{2}^{E_{22}}$ have constant signs on $%
R_{20}\cup X_{-}\cup R_{10}^{+}.$ But 
\begin{equation*}
\lambda _{2}^{E_{21}}\left( \mu _{1},0\right) =-\frac{1}{P}\mu _{1}\left(
\gamma \delta _{1}-P\right) <0
\end{equation*}
and $\lambda _{2}^{E_{22}}\left( \mu _{1},0\right) =\mu _{1}<0$ if $\mu
_{1}<0.$ Using $\lambda _{1}^{E_{21}}>0$ and $\lambda _{1}^{E_{22}}<0$
whenever $E_{21}$ and $E_{22}$ are proper, it follows that $E_{21}$ is a
saddle and $E_{22}$ an attractor on $R_{20}.$ On $X_{-}\cup R_{10}^{+},$ $%
E_{21}$ continues to remain a saddle while $E_{22}$ vanishes (it becomes a
virtual point with $\xi _{22}<0$ ).

On $T_{3},$ $E_{3}$ collides to $E_{21}$ and $\left. \lambda
_{2}^{E_{21}}\right\vert _{T_{3}}=0.$ Thus $\lambda _{2}^{E_{21}}$ changes
its sign when $\left( \mu _{1},\mu _{2}\right) $ crosses $T_{3}$ and becomes
positive on $R_{10}^{-}$ if $\mu _{1}\leq 0,$ because $\lambda
_{2}^{E_{21}}\left( 0,\mu _{2}\right) =\frac{\gamma }{P}\sqrt{-\mu _{2}P}%
\left( 1+O\left( \mu _{2}\right) \right) >0$ and $\left( 0,\mu _{2}\right)
\in R_{10}^{-}$ for $\mu _{2}<0.$ Therefore, $E_{21}$ is a repeller on $%
R_{10}^{-}$ if $\mu _{1}\leq 0.$

At this step, it is important to check if $\lambda _{2}^{E_{21}}$ changes
its sign when $\left( \mu _{1},\mu _{2}\right) $ crosses $T_{3}^{+}$ because 
$\lambda _{2}^{E_{21}}\lambda _{2}^{E_{22}}=0$ on $T_{3}^{+}.$ We observe
this does not happen because 
\begin{equation*}
\left. \lambda _{2}^{E_{21}}\right\vert _{T_{3}^{+}}=-\frac{1}{P}\mu
_{1}\left( \gamma \delta _{1}-2P\right) \left( 1+O\left( \mu _{1}\right)
\right) >0
\end{equation*}
if $\mu _{1}>0.$ Thus, $E_{21}$ is a repeller on $R_{10}^{-},$ either for $%
\mu _{1}\leq 0$ or $\mu _{1}>0.$

b) If $\delta _{1}<0,$ then $R_{20}\subset \left\{ \left( \mu _{1},\mu
_{2}\right) \in 
\mathbb{R}
^{2}\left\vert \mu _{1}>0\right. ,\mu _{2}>0\right\} .$ Proceeding as in a),
we have 
\begin{equation*}
\lambda _{2}^{E_{21}}\left( \mu _{1},0\right) =-\frac{1}{P}\mu _{1}\left(
\gamma \delta _{1}-P\right) \left( 1+O\left( \mu _{1}\right) \right) >0
\end{equation*}
and $\lambda _{2}^{E_{22}}\left( \mu _{1},0\right) =\mu _{1}>0$ if $\mu
_{1}>0,$ respectively, $\lambda _{1}^{E_{21}}>0$ and $\lambda
_{1}^{E_{22}}<0 $ whenever $E_{21}$ and $E_{22}$ exist. Thus, $E_{21}$ is a
repeller and $E_{22}$ a saddle on $R_{20}.$

$E_{22}$ vanishes on $R_{10}.$ As above, $\left. \lambda
_{2}^{E_{21}}\right\vert _{T_{3}^{+}}\neq 0,$ which implies that $\lambda
_{2}^{E_{21}}$ does not change its sign when $\left( \mu _{1},\mu
_{2}\right) $ crosses $T_{3}^{+}.$ Thus, $E_{21}$ remains a repeller on $%
X_{+}\cup R_{10}^{-}.$ Notice that $\lambda _{2}^{E_{21}}\left( 0,\mu
_{2}\right) =\frac{\gamma }{P}\sqrt{-\mu _{2}P}\left(1+O\left(\mu_{2}\right)%
\right)>0$ if $\mu _{2}<0.$

On the other hand, $\lambda _{2}^{E_{21}}$ changes its sign when $\left( \mu
_{1},\mu _{2}\right) $ crosses $T_{3}$ because $\left. \lambda
_{2}^{E_{21}}\right\vert _{T_{3}}=0.$ More exactly, $\lambda _{2}^{E_{21}}<0$
and $\lambda _{1}^{E_{21}}>0$ on $R_{10}^{+}.$ Thus, $E_{21}$ is a saddle on 
$R_{10}^{+}.$ $\square $

\bigskip

\begin{figure}[htpb]
\begin{center}
\includegraphics[width=0.3\textwidth]{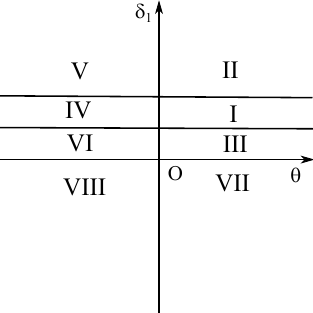}
\end{center}
\caption{When $\protect\theta\neq0$ and $\protect\delta=0,$ eight cases in
the $\protect\delta_{1}\protect\theta$- plane lead to eight bifurcation
diagrams . }
\label{fig5}
\end{figure}

\begin{figure}[htpb]
\begin{center}
\includegraphics[width=1\textwidth]{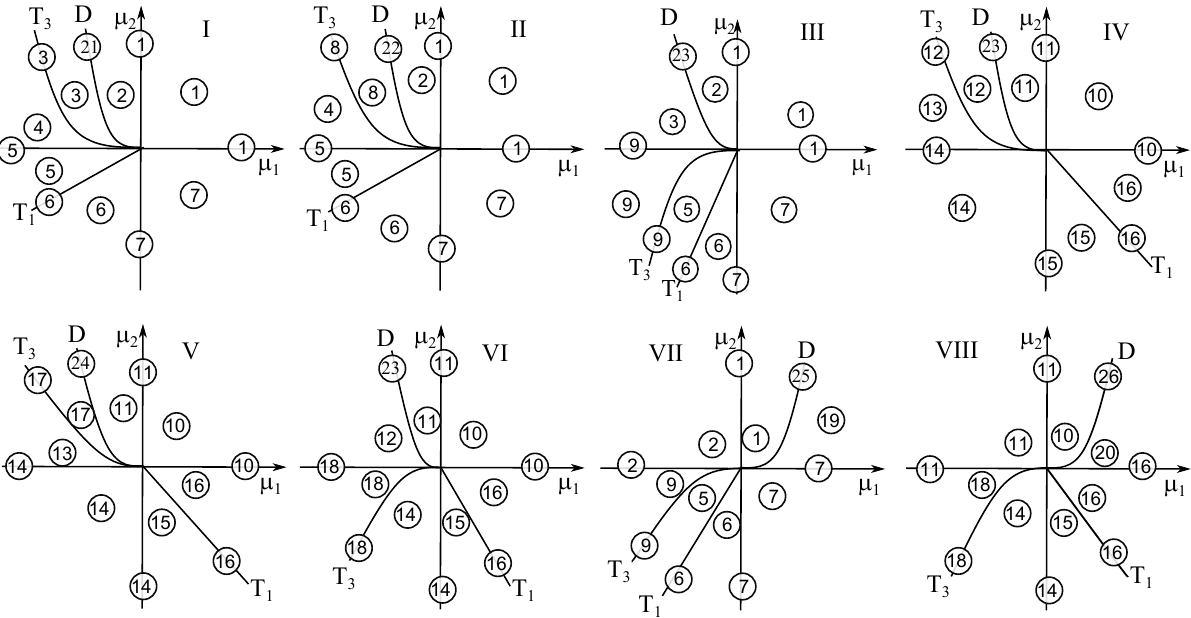}
\end{center}
\caption{Bifurcation diagrams for $\protect\theta\neq0$ and $\protect\delta=0,$ corresponding to the eight cases I-VIII.}
\label{fig6}
\end{figure}

\begin{figure}[tbph]
\begin{center}
\includegraphics[width=1\textwidth]{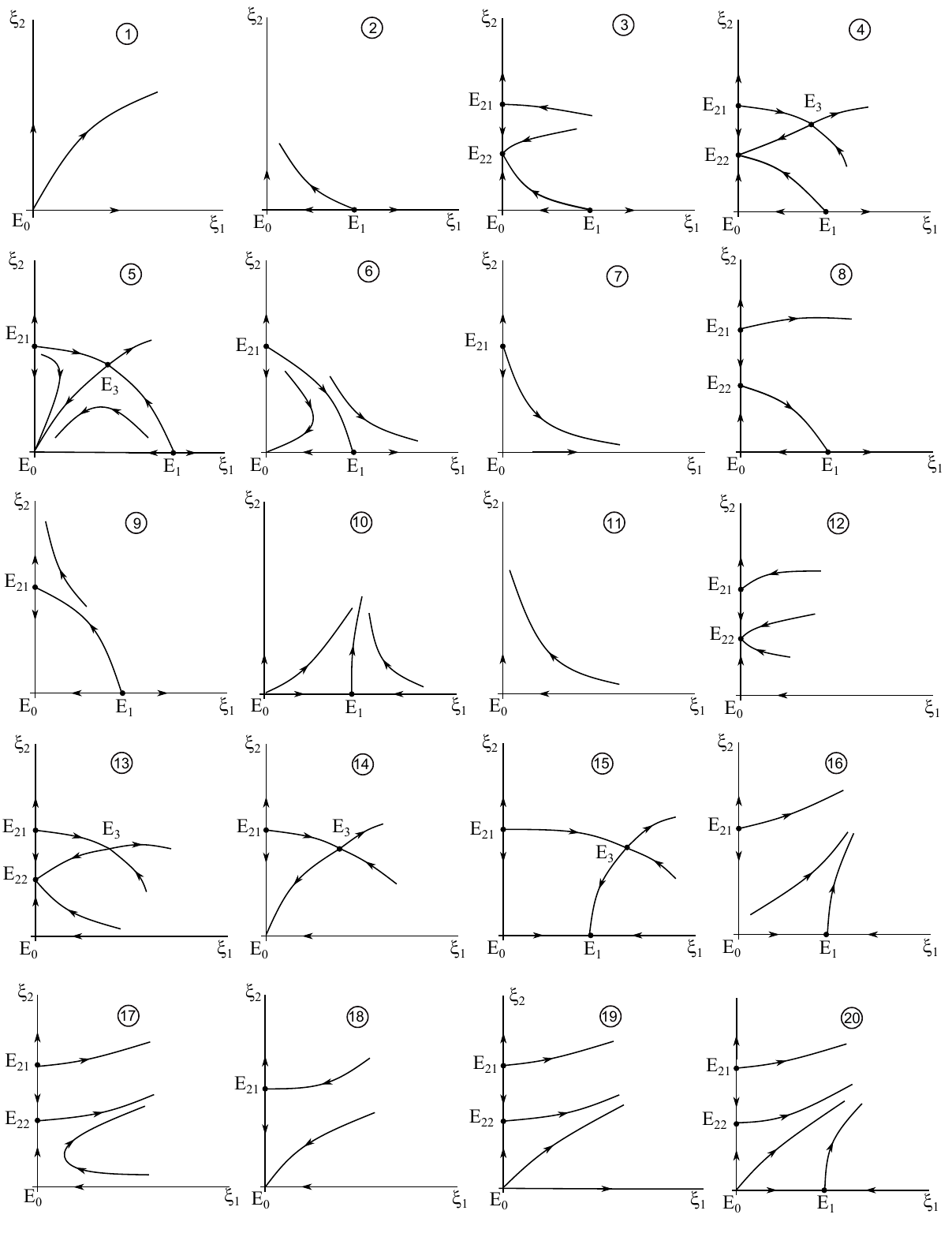}
\end{center}
\caption{Phase portraits corresponding to the bifurcation diagrams I-VIII,
when $\protect\theta\neq 0$ and $\protect\delta=0.$ }
\label{fig7}
\end{figure}

\begin{figure}[tbph]
\begin{center}
\includegraphics[width=12 cm,height=8 cm]{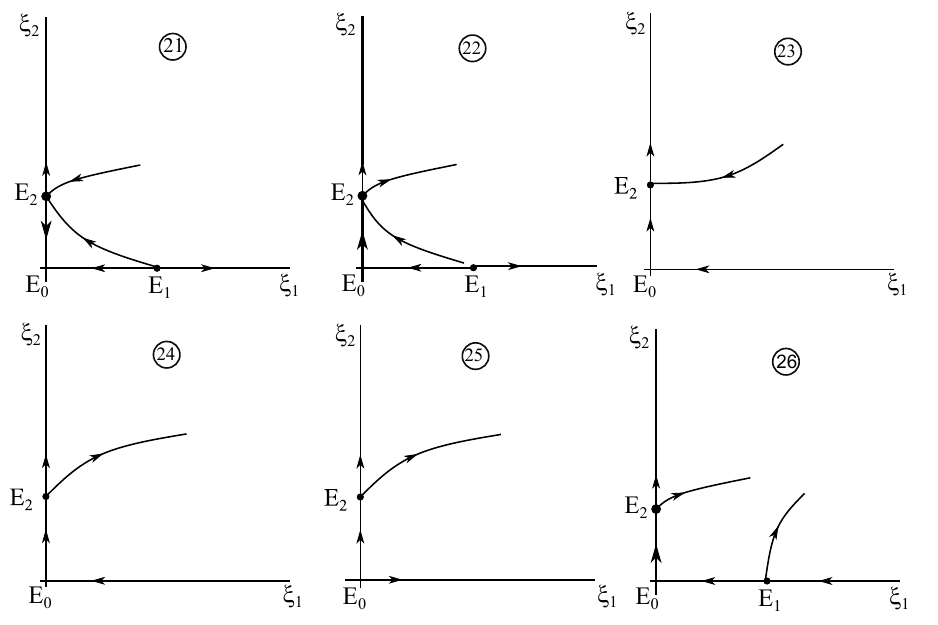}
\end{center}
\caption{Phase portraits on the saddle-node bifurcation curve, when $\protect\theta\neq0$ and $\protect\delta=0.$ }
\label{fig8}
\end{figure}

From the above results, eight different cases arise (Fig.\ref{fig5}) in
terms of $\delta _{1}$ and $\theta ,$ each one leading to a bifurcation
diagram. There are 20 different regions in the eight bifurcation diagrams
(Fig.\ref{fig6}). In Table \ref{tabel3} we summarized the type of each
equilibrium point from the 20 regions. The phase portraits corresponding to
the 20 regions are depicted in Fig.\ref{fig7}. The phase portraits on the
saddle node curve are presented in Fig.\ref{fig8}.

\bigskip

I: $\theta >0,$ $\delta _{1}>0,$ $\gamma \delta _{1}-P>0$ and $\gamma \delta
_{1}-2P<0,$ Fig.\ref{fig5} (I)

II: $\theta >0,$ $\delta _{1}>0,$ $\gamma \delta _{1}-P>0$ and $\gamma
\delta _{1}-2P>0,$ Fig.\ref{fig5} (II)

III: $\theta >0,$ $\delta _{1}>0,$ $\gamma \delta _{1}-P<0$ and $\gamma
\delta _{1}-2P<0,$ Fig.\ref{fig5} (III)

IV: $\theta <0,$ $\delta _{1}>0,$ $\gamma \delta _{1}-P>0$ and $\gamma
\delta _{1}-2P<0,$ Fig.\ref{fig5} (IV)

V: $\theta <0,$ $\delta _{1}>0,$ $\gamma \delta _{1}-P>0$ and $\gamma \delta
_{1}-2P>0,$ Fig.\ref{fig5} (V)

VI: $\theta <0,$ $\delta _{1}>0,$ $\gamma \delta _{1}-P<0$ and $\gamma
\delta _{1}-2P<0,$ Fig.\ref{fig5} (VI)

VII: $\theta >0,$ $\delta _{1}<0,$ Fig.\ref{fig5} (VII)

VIII: $\theta <0,$ $\delta _{1}<0,$ Fig.\ref{fig5} (VIII) 

\begin{table}[htpb]
\centering
\begin{tabular}{l|llllllllllllllllllll}
& $1$ & $2$ & $3$ & $4$ & $5$ & $6$ & $7$ & $8$ & $9$ & $10$ & $11$ & $12$ & 
$13$ & $14$ & $15$ & $16$ & $17$ & $18$ & $19$ & $20$ \\ \hline
$O$ & $r$ & $s$ & $s$ & $s$ & $a$ & $a$ & $s$ & $s$ & $a$ & $r$ & $s$ & $s$
& $s$ & $a$ & $s$ & $s$ & $s$ & $a$ & $r$ & $r$ \\ 
$E_{1}$ & $-$ & $r$ & $r$ & $r$ & $r$ & $s$ & $-$ & $r$ & $r$ & $s$ & $-$ & $%
-$ & $-$ & $-$ & $a$ & $s$ & $-$ & $-$ & $-$ & $s$ \\ 
$E_{21}$ & $-$ & $-$ & $s$ & $r$ & $r$ & $r$ & $r$ & $r$ & $s$ & $-$ & $-$ & 
$s$ & $r$ & $r$ & $r$ & $r$ & $r$ & $s$ & $r$ & $r$ \\ 
$E_{22}$ & $-$ & $-$ & $a$ & $a$ & $-$ & $-$ & $-$ & $s$ & $-$ & $-$ & $-$ & 
$a$ & $a$ & $-$ & $-$ & $-$ & $s$ & $-$ & $s$ & $s$ \\ 
$E_{3}$ & $-$ & $-$ & $-$ & $s$ & $s$ & $-$ & $-$ & $-$ & $-$ & $-$ & $-$ & $%
-$ & $s$ & $s$ & $s$ & $-$ & $-$ & $-$ & $-$ & $-$%
\end{tabular}%
\caption{\textit{The types of the equilibrium points of system \eqref{gl8}
for $\protect\theta \neq 0$ and $\protect\delta =0$ on different regions of
the bifurcation diagrams I-VIII; the abbreviations s, a, r stand for saddle,
attractor, repeller, respectively.}}
\label{tabel3}
\end{table}

\newpage 

\section{The degeneracy of the system at $\protect\theta \left( 0\right) =0$
and $\protect\delta \left( 0\right) \neq 0$}

\bigskip Since $\theta \left( 0\right) =\theta =0,$ we have 
\begin{equation*}
\theta \left( \mu \right) =\theta _{1}\mu _{1}+\theta _{2}\mu _{2}+O\left(
\left\vert \mu \right\vert ^{2}\right)
\end{equation*}%
and similarly for the other coefficients; $\gamma \left( 0\right) =\gamma
>0. $ Assume $\frac{\partial \theta \left( 0\right) }{\partial \mu _{1}}%
\overset{not}{=}\theta _{1}\neq 0,$ $\frac{\partial \theta \left( 0\right) }{%
\partial \mu _{2}}\overset{not}{=}\theta _{2}\neq 0,$ $N\left( 0\right) 
\overset{not}{=}N\neq 0,$ respectively, $\delta =\delta \left( 0\right) \neq
0.$

The system (\ref{gl8}) has two equilibrium points $E_{11}\left( \xi
_{11},0\right) $ and $E_{12}\left( \xi _{12},0\right) $ lying on the $\xi
_{1}-$axis, given by

\begin{equation*}
\xi _{11}=-\frac{1}{2N\left( \mu \right) }\left( \theta \left( \mu \right) -%
\sqrt{\Delta ^{\prime }\left( \mu \right) }\right) \text{ and }\xi _{12}=-%
\frac{1}{2N\left( \mu \right) }\left( \theta \left( \mu \right) +\sqrt{%
\Delta ^{\prime }\left( \mu \right) }\right) ,
\end{equation*}%
where $\Delta ^{\prime }\left( \mu \right) =\theta ^{2}\left( \mu \right) -$ 
$4\mu _{1}N\left( \mu \right) \geq 0.$

\begin{remark}
\label{rtheta1} The eigenvalues of $E_{0}$ are $\mu _{1}$ and $\mu _{2},$
while of $E_{2}$ they are $(\mu _{1}-\frac{\gamma }{\delta }\mu _{2})\left(
1+O\left( \left\vert \mu \right\vert \right) \right) $ and $-\mu _{2}.$
\end{remark}

For $\left\vert \mu \right\vert $ sufficiently small, denote by 
\begin{equation*}
D=\left\{ \left( \mu _{1},\mu _{2}\right) \in 
\mathbb{R}
^{2}\left\vert \mu _{1}=\frac{\theta _{2}^{2}}{4N}\mu _{2}^{2}\left(
1+O\left( \mu _{2}\right) \right) \right. \right\}
\end{equation*}%
the bifurcation curve $\Delta ^{\prime }\left( \mu \right) =0.$ Similar to
the previous degenerate case, $E_{11}\left( \xi _{11},0\right) $ and $%
E_{12}\left( \xi _{12},0\right) $ come into existence by a saddle-node
bifurcation, which is proved in the next theorem.

\begin{theorem}
If $\theta _{1}\theta _{2}\delta N\left( 2N\gamma -\theta _{2}\right) \neq 0$
and $\gamma >0,$ then $D\cap \left\{ (\mu _{1},\mu _{2})\mid \mu
_{2}>0\right\} $ and $D\cap \left\{ (\mu _{1},\mu _{2})\mid \mu
_{2}<0\right\} $ are two saddle-node bifurcation curves.
\end{theorem}

\textit{Proof.} Consider first the branch $D^{+}=D\cap \left\{ \mu
_{2}>0\right\} .$ Write the system in the form $\frac{d\xi }{dt}=f\left( \xi
,\mu \right) ,$ with $\xi =\left( \xi _{1},\xi _{2}\right) ,$ $f=\left(
f_{1},f_{2}\right) $ and $\mu =\left( \mu _{1},\mu _{2}\right) .$ The
equilibrium $\left( \xi _{1},0\right) $ satisfies $\mu _{1}+\theta \left(
\mu \right) \xi _{1}+N\left( \mu \right) \xi _{1}^{2}=0.$ Then $f\left( \xi
_{0},\mu _{0}\right) =\left( 0,0\right) $ for $\xi _{0}=\left( \xi
_{1},0\right) $ and $\mu _{0}=\left( \mu _{1},\mu _{2}\right) \in D^{+}.$
Notice that, $\xi _{11}=\xi _{12}=\xi _{1}=-\frac{\theta \left( \mu \right) 
}{2N\left( \mu \right) }$ whenever $\mu _{0}\in D^{+}.$

Assume that $\mu _{2}>0$ is fixed while $\mu _{1}$ varies, thus, $\mu _{1}$
is the bifurcation parameter. The Jacobian matrices 
\begin{equation*}
A=Df\left( \xi _{0},\mu _{0}\right) =\left( 
\begin{array}{cc}
0 & \xi _{1}\left( \gamma \left( \mu _{0}\right) +M\left( \mu _{0}\right)
\xi _{1}\right) \\ 
0 & \frac{1}{\gamma \left( \mu _{0}\right) }\left( \xi _{1}+\gamma \left(
\mu _{0}\right) \mu
_{2}+R\left(\mu_{0}\right)\gamma\left(\mu_{0}\right)\xi_{1}^{2}\right)%
\end{array}%
\right)
\end{equation*}%
and $A^{T}$ have both an eigenvalue $\lambda =0,$ with the corresponding
eigenvectors $v$ for $A$ and $w$ for $A^{T},$ where $v=\left( 
\begin{array}{c}
1 \\ 
0%
\end{array}%
\right) \ $and $w=\left( 
\begin{array}{c}
-\frac{\xi _{1}+\gamma \left( \mu _{0}\right) \mu
_{2}+R\left(\mu_{0}\right)\gamma\left(\mu_{0}\right)\xi_{1}^{2}}{\gamma
\left( \mu _{0}\right) \xi _{1}\left( \gamma \left( \mu _{0}\right) +M\left(
\mu _{0}\right) \xi _{1}\right) } \\ 
1%
\end{array}%
\right) .$ We have also $f_{\mu _{1}}=\left( 
\begin{array}{cc}
\frac{\partial f_{1}}{\partial \mu _{1}} & \frac{\partial f_{2}}{\partial
\mu _{1}}%
\end{array}%
\right) ^{T}=\left( 
\begin{array}{cc}
\xi _{1}\left( 1+\frac{\partial \theta \left( \mu \right) }{\partial \mu _{1}%
}\xi _{1}+\frac{\partial N\left( \mu \right) }{\partial \mu _{1}}\xi
_{1}^{2}\right) & 0%
\end{array}%
\right) ^{T}$ at $\left( \xi _{1},0\right) .$ These lead to

\begin{equation*}
w^{T}f_{\mu _{1}}\left( \xi _{0},\mu _{0}\right) =-\mu _{2}\frac{2N\gamma 
\text{ }-\text{ }\theta _{2}}{2N\gamma ^{2}}\left( 1+O\left( \mu _{2}\right)
\right) \neq 0,
\end{equation*}%
respectively,%
\begin{equation*}
w^{T}\left[ D^{2}f\left( \xi _{0},\mu _{0}\right) \left( v,v\right) \right]
=-\mu _{2}\frac{2N\gamma -\theta _{2}}{\gamma ^{2}}\left( 1+O\left( \mu
_{2}\right) \right) \neq 0,
\end{equation*}%
which confirm the claim. Notice that $w=\allowbreak \left( 
\begin{array}{cc}
-\frac{1}{\gamma ^{2}\theta _{2}}\left( \theta _{2}-2N\gamma \right) \left(
1+O\left( \mu _{2}\right) \right) & 1%
\end{array}%
\right) ^{T}$ and $f_{\mu _{1}}\left( \xi _{0},\mu _{0}\right) =\left( 
\begin{array}{cc}
-\frac{\theta _{2}\mu _{2}}{2N}\left( 1+O\left( \mu _{2}\right) \right) & 0%
\end{array}%
\right) ^{T}$ at $\left( \xi _{0},\mu _{0}\right) .$ For $D\cap \left\{ \mu
_{2}<0\right\} $ the proof is similar. $\square $

\begin{theorem}
\label{ttheta1} If $\theta _{1}\theta _{2}\delta N\left( \theta _{2}-N\gamma
\right) \neq 0$ and $\gamma >0,$ the system (\ref{gl8}) has an equilibrium
point of the form $E_{3}\left( \xi _{1},\xi _{2}\right) ,$ which is a saddle
whenever it exists.
\end{theorem}

\textit{Proof.} From the Implicit Functions Theorem, the system (\ref{ec2})
has a unique solution of the form

\begin{equation}
\xi _{1}=\left( \delta \mu _{1}-\gamma \mu _{2}\right) (1+O\left( |\mu
|\right) )\text{ and }\xi _{2}=\left( -\text{ }\frac{\mu _{1}}{\gamma }\text{
}+\left( \theta _{2}\text{ }-\text{ }N\gamma \right) \mu _{2}^{2}\right)
(1+O\left( \left\vert \mu \right\vert \right) ),  \label{eqtheta2}
\end{equation}%
for all $\left\vert \mu \right\vert $ sufficiently small. Thus, $E_{3}\left(
\xi _{1},\xi _{2}\right) $ with $\xi _{1,2}$ given by (\ref{eqtheta2}) is an
equilibrium point of (\ref{gl8}). It bifurcates from $O$ along the curves 
\begin{equation*}
T_{2}=\left\{ \left( \mu _{1},\text{ }\mu _{2}\right) \in 
\mathbb{R}
^{2}\left\vert \text{ }\mu _{2}=\frac{\delta }{\gamma }\mu _{1}\left(
1+O\left( \mu _{1}\right) \right) \right. ,\delta \mu _{2}<0\right\} ,
\end{equation*}
respectively, 
\begin{equation}
T_{4}=\left\{ \left( \mu _{1},\mu _{2}\right) \in 
\mathbb{R}
^{2}\left\vert \mu _{1}=\gamma \left( \theta _{2}-N\gamma \right) \mu
_{2}^{2}\left( 1+O\left( \mu _{2}\right) \right) \right. ,\mu _{2}<0\right\}
.  \label{t4}
\end{equation}

$E_{3}$ collides to $E_{2}\left( 0,-\frac{1}{\delta }\mu _{2}\left(
1+O\left( \left\vert \mu \right\vert \right) \right) \right) $ on $T_{2},$
respectively, $E_{11}\left( \xi _{11},0\right) $ or $E_{12}\left( \xi
_{12},0\right) $ on $T_{4}.$ We call $E_{3}$ \textit{trivial} in these
cases, otherwise nontrivial. In its lowest terms, $E_{3}$ reads 
\begin{equation*}
E_{3}\left( \delta \mu _{1}-\gamma \mu _{2},-\frac{\mu _{1}}{\gamma }+\left(
\theta _{2}-N\gamma \right) \mu _{2}^{2}\right) .
\end{equation*}

Thus, $E_{3}$ is a sadlle whenever it exists in $Q_{1}$ and is non-trivial,
since its eigenvalues also satisfy (\ref{ell}). $\square $

\begin{theorem}
\label{ttheta2} Assume $\delta \theta _{2}N\neq 0,$ $\gamma >0$ and $\left(
\mu _{1},\mu _{2}\right) \in T_{4}$ with $\left\vert \mu \right\vert $
sufficiently small. Then,

a) if $\theta _{2}-2N\gamma <0,$ $E_{3}$ collides to $E_{11}\left( -\gamma
\mu _{2},0\right) ,$ $\left. \lambda _{1}^{E_{11}}\right\vert _{T_{4}}>0$
and $\left. \lambda _{2}^{E_{11}}\right\vert _{T_{4}}=0.$ Moreover, $%
E_{12}\left( \mu _{2}\frac{N\gamma -\theta _{2}}{N},0\right) $ is an
attractor, whenever it is proper.

b) if $\theta _{2}-2N\gamma >0,$ $E_{3}$ collides to $E_{12}\left( -\gamma
\mu _{2},0\right) ,$ $\left. \lambda _{1}^{E_{12}}\right\vert _{T_{4}}<0$
and $\left. \lambda _{2}^{E_{12}}\right\vert _{T_{4}}=0.$ Moreover, $%
E_{11}\left( \mu _{2}\frac{N\gamma -\theta _{2}}{N},0\right) $ is a
repeller, whenever it is proper.
\end{theorem}

\textit{Proof.} The eigenvalues of $E_{11}$ are $\lambda _{1}^{E_{11}}=-2\mu
_{1}-\mu _{2}\theta _{2}\xi _{11}$ and $\lambda _{2}^{E_{11}}=\mu _{2}+\frac{%
1}{\gamma }\xi _{11}+R\xi_{11}^{2},$ while of $E_{12}$ they are $\lambda
_{1}^{E_{22}}=-2\mu _{1}-\mu _{2}\theta _{2}\xi _{12}$ and $\lambda
_{2}^{E_{22}}=\mu _{2}+\frac{1}{\gamma }\xi _{12}+R\xi_{12}^{2}.$ These lead
to

\begin{equation}
\lambda _{1}^{E_{11}}=\xi _{11}\sqrt{\Delta ^{\prime }\left( \mu \right) }>0%
\text{ and }\lambda _{1}^{E_{12}}=-\xi _{12}\sqrt{\Delta ^{\prime }\left(
\mu \right) }<0.  \label{ll11}
\end{equation}

Since it is difficult to study the signs of $\lambda _{2}^{E_{11}}$ and $%
\lambda _{2}^{E_{12}}$ as separate terms, we link them together through the
product 
\begin{equation}
\lambda _{2}^{E_{11}}\lambda _{2}^{E_{12}}=\frac{\mu _{1}-\gamma \left(
\theta _{2}-N\gamma \right) \mu _{2}^{2}}{\gamma ^{2}N}\left( 1+O\left( \mu
\right) \right) .  \label{c4}
\end{equation}

a) Assume $\theta _{2}-2N\gamma <0.$ Then, $E_{3}$ collides to $E_{11}\left(
-\gamma \mu _{2},0\right) $ for $\left\vert \mu \right\vert $ sufficiently
small and $\left( \mu _{1},\mu _{2}\right) \in T_{4}.$

So, $\lambda _{2}^{E_{11}}=0$ and $\lambda _{1}^{E_{11}}=\xi _{11}\sqrt{%
\Delta ^{\prime }\left( \mu \right) }>0,$ on $T_{4}.$ In addition, $%
E_{12}\left( \mu _{2}\frac{N\gamma -\theta _{2}}{N},0\right) $ has its
eigenvalues $\left. \lambda _{1}^{E_{12}}\right\vert _{T_{4}}=-\xi _{12}%
\sqrt{\Delta ^{\prime }}<0$ and 
\begin{equation*}
\left. \lambda _{2}^{E_{12}}\right\vert _{T_{4}}=-\mu _{2}\frac{\theta
_{2}-2N\gamma }{N}\left( 1+O\left( \mu _{2}\right) \right) <0,
\end{equation*}%
since $\mu _{2}<0$ on $T_{4}.$ Thus, $E_{12}$ is an attractor on $T_{4}.$
Since $\left. \lambda _{2}^{E_{12}}\right\vert _{T_{4}}\neq 0$ and $\left.
\lambda _{2}^{E_{11}}\right\vert _{T_{4}}=0,$ the curve $\lambda
_{2}^{E_{11}}\lambda _{2}^{E_{12}}=0$ given by (\ref{c4}) coincides to $%
T_{4}.$

b) Assume $\theta _{2}-2N\gamma >0.$ Then, $E_{3}$ collides to $E_{12}\left(
-\gamma \mu _{2},0\right) \ $for $\left\vert \mu \right\vert $ sufficiently
small and $\left( \mu _{1},\mu _{2}\right) \in T_{4}.$ Also, $\lambda
_{1}^{E_{12}}\lambda _{2}^{E_{12}}=0$ on $T_{4}$ by (\ref{ell}), which
yields $\lambda _{2}^{E_{12}}=0$ on $T_{4},$ because $\lambda
_{1}^{E_{12}}\neq 0$ on $T_{4}.$ In addition, $E_{11}\left( \mu _{2}\frac{%
N\gamma -\theta _{2}}{N},0\right) $ has its eigenvalues $\left. \lambda
_{1}^{E_{11}}\right\vert _{T_{4}}=\xi _{11}\sqrt{\Delta ^{\prime }\left( \mu
\right) }>0$ and 
\begin{equation*}
\left. \lambda _{2}^{E_{11}}\right\vert _{T_{4}}=-\mu _{2}\frac{\theta
_{2}-2N\gamma }{N}\left( 1+O\left( \mu _{2}\right) \right) >0,
\end{equation*}%
since $\mu _{2}<0$ on $T_{4}.$ Thus, $E_{11}$ is a repeller on $T_{4}.$
Since $\lambda _{2}^{E_{11}}\lambda _{2}^{E_{12}}=0$ on $T_{4}$ and $\left.
\lambda _{2}^{E_{11}}\right\vert _{T_{4}}\neq 0,$ it follows again that the
curve $\lambda _{2}^{E_{11}}\lambda _{2}^{E_{12}}=0$ coincides to $T_{4}.$ $%
\square $

Define the following regions 
\begin{equation*}
R_{00}=\left\{ \left( \mu _{1},\mu _{2}\right) \in 
\mathbb{R}
^{2}\left\vert \Delta ^{\prime }\left( \mu \right) >0,\mu _{1}>0,\text{ }%
\theta _{2}\mu _{2}>0\right. \right\} \cup \left\{ \left( \mu _{1},\mu
_{2}\right) \in 
\mathbb{R}
^{2}\left\vert \Delta ^{\prime }\left( \mu \right) <0\right. \right\} ,
\end{equation*}

\begin{equation*}
R_{10}=\left\{ \left( \mu _{1},\mu _{2}\right) \in 
\mathbb{R}
^{2}\left\vert \Delta ^{\prime }\left( \mu \right) >0,\mu _{1}<0\right.
\right\} \text{ and }R_{20}=\left\{ \left( \mu _{1},\mu _{2}\right) \in 
\mathbb{R}
^{2}\left\vert \Delta ^{\prime }\left( \mu \right) >0,\mu _{1}>0,\text{ }%
\theta _{2}\mu _{2}<0\right. \right\} .
\end{equation*}

Both points $E_{11}$ and $E_{22}$ are proper in the region $R_{20},$because $%
\xi _{11}>\xi _{12}>0,$ while $R_{10}$ contains only $E_{11}$ proper, since $%
\xi _{11}>0>\xi _{12}.$ On $R_{00},$ none of the two equilibria survive as
proper points (they become virtual) because either $\xi _{11}<0$ and $\xi
_{12}<0$ or $\xi _{11}$ and $\xi _{12}$ are not real numbers.

Whenever $\theta _{2}-N\gamma >0,$ denote by 
\begin{equation*}
R_{20}^{-}=R_{20}\cap \left\{ \mu _{1}<\gamma \left( \theta _{2}-N\gamma
\right) \mu _{2}^{2}\right\} \text{ and }R_{20}^{+}=R_{20}\cap \left\{ \mu
_{1}>\gamma \left( \theta _{2}-N\gamma \right) \mu _{2}^{2}\right\} ,
\end{equation*}
the regions from $R_{20}$ to the left, respectively, the right of $T_{4}.$
Notice that 
\begin{equation*}
R_{20}=R_{20}^{-}\cup T_{4}\cup R_{20}^{+}.
\end{equation*}

\begin{theorem}
\label{ttheta3} Assume $\delta \theta _{2}\neq 0,$ $\gamma >0,$ $N>0$ and $%
\theta _{2}-N\gamma >0.$ Then,

a) if $\theta _{2}-2N\gamma <0,$ then $E_{11}$ is a repeller on $R_{10}\cup
Y_{-}\cup R_{20}^{-}$ and a saddle on $R_{20}^{+},$ while $E_{12}$ is an
attractor on $R_{20}.$

b) if $\theta _{2}-2N\gamma >0,$ then $E_{11}$ is a repeller on $R_{10}\cup
Y_{-}\cup R_{20},$ while $E_{12}$ is a saddle on $R_{20}^{+}$ and an
attractor on $R_{20}^{-}.$
\end{theorem}

\textit{Proof.} a) Notice first that $\theta _{2}>0$ and $0<N\gamma <\theta
_{2}<2N\gamma ,$ thus, 
\begin{equation*}
R_{20}\subset \left\{ \left( \mu _{1},\mu _{2}\right) \in 
\mathbb{R}
^{2}\left\vert \mu _{1}>0\right. ,\mu _{2}<0\right\} .
\end{equation*}%
From Theorem \ref{ttheta2}, $E_{3}$ collides to $E_{11}\left( -\gamma \mu
_{2},0\right) $ on $T_{4}$ and $\left. \lambda _{2}^{E_{11}}\right\vert
_{T_{4}}=0.$ Thus, $\lambda _{2}^{E_{11}}\left( \mu _{1},\mu _{2}\right) $
keeps constant sign outside $T_{4}$ and changes its sign when $\left( \mu
_{1},\mu _{2}\right) $ crosses $T_{4}.$ Since 
\begin{equation*}
\lambda _{2}^{E_{11}}\left( \mu _{1},0\right) =\frac{1}{N\gamma }\sqrt{-\mu
_{1}N}\left( 1+O\left( \mu _{1}\right) \right) >0
\end{equation*}%
if $\mu _{1}<0,$ it follows that $\lambda _{2}^{E_{11}}>0$ on $R_{20}^{-}$
and $\lambda _{2}^{E_{11}}<0$ on $R_{20}^{+},$ while $\lambda
_{1}^{E_{11}}=\xi _{11}\sqrt{\Delta ^{\prime }\left( \mu \right) }>0$ on $%
R_{20}.$ Thus, $E_{11}$ is a saddle on $R_{20}^{+}$ and a repeller on $%
R_{20}^{-}.$ It remains a repeller on $R_{10}\cup X_{-}$ based on the same
reasons.

On the other hand, since $\left. \lambda _{2}^{E_{12}}\right\vert
_{T_{4}}\neq 0,$ $\lambda _{2}^{E_{12}}$ keeps constant sign $R_{20},$ which
is negative by Theorem \ref{ttheta2}. From $\lambda _{1}^{E_{12}}=-\xi _{12}%
\sqrt{\Delta ^{\prime }\left( \mu \right) }<0,$ it follows that $E_{12}$ is
an attractor on $R_{20}.$

b) We have in this case $\theta _{2}>0$ and $0<2N\gamma <\theta _{2},$ which
yield 
\begin{equation*}
R_{20}\subset \left\{ \left( \mu _{1},\mu _{2}\right) \in 
\mathbb{R}
^{2}\left\vert \mu _{1}>0\right. ,\mu _{2}<0\right\} .
\end{equation*}%
From Theorem \ref{ttheta2}, $E_{3}$ collides to $E_{12}\left( -\gamma \mu
_{2},0\right) $ on $T_{4}$ and $\left. \lambda _{2}^{E_{12}}\right\vert
_{T_{4}}=0,$ thus, $\lambda _{2}^{E_{12}}\left( \mu _{1},\mu _{2}\right) $
changes its sign when $\left( \mu _{1},\mu _{2}\right) $ crosses $T_{4}.$
More exactly, $\lambda _{2}^{E_{12}}>0$ on $R_{20}^{+}$ because 
\begin{equation*}
\left. \lambda _{2}^{E_{12}}\right\vert _{D}=-\frac{1}{2N\gamma }\mu
_{2}\left( \theta _{2}-2N\gamma \right) \left( 1+O\left( \mu _{2}\right)
\right) >0,
\end{equation*}%
which, in turn, yields $\lambda _{2}^{E_{12}}<0$ on $R_{20}^{-}.$ Thus, $%
E_{12}$ is a saddle on $R_{20}^{+}$ and an attractor on $R_{20}^{-},$
because $\lambda _{1}^{E_{12}}<0$ on $R_{20}.$

Further, $\lambda _{2}^{E_{11}}\neq 0$ on $R_{10}\cup Y_{-}\cup R_{20}$ and $%
\left. \lambda _{2}^{E_{11}}\right\vert _{T_{4}}=-\mu _{2}\frac{\theta
_{2}-2N\gamma }{N}\left(1+O\left(\mu_{2}\right)\right)>0,$ yield $\lambda
_{2}^{E_{11}}>0$ on $R_{10}\cup Y_{-}\cup R_{20}.$ Thus, $E_{11}$ is a
repeller on $R_{10}\cup Y_{-}\cup R_{20}.$ $\square $

\bigskip

Whenever $\theta _{2}-N\gamma <0,$ denote by 
\begin{equation*}
R_{10}^{+}=R_{10}\cap \left\{ \mu _{1}>\gamma \left( \theta _{2}-N\gamma
\right) \mu _{2}^{2}\right\} \text{ and }R_{10}^{-}=R_{10}\cap \left\{ \mu
_{1}<\gamma \left( \theta _{2}-N\gamma \right) \mu _{2}^{2}\right\} ,
\end{equation*}
such that $R_{10}=R_{10}^{+}\cup T_{4}\cup R_{10}^{-}.$ Denote by 
\begin{equation*}
T_{4}^{+}=\left\{ \left( \mu _{1},\mu _{2}\right) \in 
\mathbb{R}
^{2}\left\vert \mu _{1}=\gamma \left( \theta _{2}-N\gamma \right) \mu
_{2}^{2}\left( 1+O\left( \mu _{1}\right) \right) \right. ,\mu _{2}>0\right\}
.
\end{equation*}

\begin{theorem}
\label{ttheta4} Assume $\delta \theta _{2}\neq 0,$ $\gamma >0,$ $N>0$ and $%
\theta _{2}-N\gamma <0.$ Then,

a) if $\theta _{2}>0,$ $E_{11}$ is a saddle and $E_{12}$ an attractor on $%
R_{20}.$ Moreover, $E_{11}$ is a saddle on $Y_{-}\cup R_{10}^{+}$ and a
repeller on $R_{10}^{-}.$

b) if $\theta _{2}<0,$ $E_{11}$ is a repeller and $E_{12}$ a saddle on $%
R_{20}.$ Moreover, $E_{11}$ is a repeller on $Y_{+}\cup R_{10}^{-}$ and a
saddle on $R_{10}^{+}.$
\end{theorem}

\textit{Proof. }a) Assume first $\theta _{2}>0,$ thus, 
\begin{equation*}
R_{20}\subset \left\{ \left( \mu _{1},\mu _{2}\right) \left\vert \mu
_{1}>0\right. ,\mu _{2}<0\right\}
\end{equation*}%
and $T_{4}\subset R_{10}.$ If $\mu _{1}>0,$ then $\lambda
_{2}^{E_{11}}\lambda _{2}^{E_{12}}=0$ only on $T_{4}\nsubseteq R_{20}\cup
Y_{-}\cup R_{10}^{+}.$ Thus, $\lambda _{2}^{E_{11}}$ and $\lambda
_{2}^{E_{12}}$ have constant signs on $R_{20}\cup Y_{-}\cup R_{10}^{+}.$ But 
\begin{equation*}
\lambda _{2}^{E_{11}}\left( 0,\text{ }\mu _{2}\right) =-\frac{1}{N\gamma }%
\mu _{2}\left( \theta _{2}-N\gamma \right) \left( 1+O\left( \mu _{2}\right)
\right) <0
\end{equation*}%
and $\lambda _{2}^{E_{12}}\left( 0,\text{ }\mu _{2}\right) =\mu _{2}<0$ if $%
\mu _{2}<0.$ Using $\lambda _{1}^{E_{11}}>0$ and $\lambda _{1}^{E_{12}}<0$
whenever $E_{11}$ and $E_{12}$ exist in $Q_{1},$ it follows that $E_{11}$ is
a saddle and $E_{12}$ an attractor on $R_{20}.$ On $Y_{-}\cup R_{10}^{+},$ $%
E_{11}$ continues to remain a saddle while $E_{12}$ vanishes (it becomes a
virtual point with $\xi _{12}<0$ ).

On $T_{4},$ $E_{3}$ collides to $E_{11}$ and $\left. \lambda
_{2}^{E_{11}}\right\vert _{T_{3}}=0.$ Thus $\lambda _{2}^{E_{11}}$ changes
its sign when $\left( \mu _{1},\mu _{2}\right) $ crosses $T_{4}$ and becomes
positive on $R_{10}^{-}$ if $\mu _{2}\leq 0,$ because 
\begin{equation*}
\lambda _{2}^{E_{11}}\left( \mu _{1},0\right) =\frac{1}{N\gamma }\sqrt{-\mu
_{1}N}\left( 1+O\left( \mu _{1}\right) \right) >0
\end{equation*}%
and $\left( \mu _{1},0\right) \in R_{10}^{-}$ for $\mu _{1}<0.$ Therefore, $%
E_{11}$ is a repeller on $R_{10}^{-}$ if $\mu _{1}\leq 0.$

Notice that $\lambda _{2}^{E_{11}}$ does not change its sign when $\left(
\mu _{1},\mu _{2}\right) $ crosses $T_{4}^{+}$ because 
\begin{equation*}
\left. \lambda _{2}^{E_{11}}\right\vert _{T_{4}^{+}}=-\frac{1}{N\gamma }\mu
_{2}\left( \theta _{2}-N\gamma \right) \left( 1+O\left( \mu _{2}\right)
\right) >0
\end{equation*}%
if $\mu _{2}<0.$ Thus, $E_{11}$ is a repeller on $R_{10}^{-},$ either for $%
\mu _{2}\leq 0$ or $\mu _{2}>0.$

b) If $\theta _{2}<0,$ then 
\begin{equation*}
R_{20}\subset \left\{ \mu _{1}>0,\mu _{2}>0\right\} .
\end{equation*}%
\newline
Proceeding as in a), we have $\lambda _{2}^{E_{11}}\left( 0,\text{ }\mu
_{2}\right) =-\frac{1}{N\gamma }\mu _{2}\left( \theta _{2}-N\gamma \right)
\left( 1+O\left( \mu _{2}\right) \right) >0$ and $\lambda
_{2}^{E_{12}}\left( 0,\text{ }\mu _{2}\right) =\mu _{2}>0$ if $\mu _{2}>0,$
respectively, $\lambda _{1}^{E_{11}}>0$ and $\lambda _{1}^{E_{12}}<0$
whenever $E_{11}$ and $E_{12}$ are proper. Thus, $E_{11}$ is a repeller and $%
E_{12}$ a saddle on $R_{20}.$

$E_{12}$ vanishes on $R_{10}.$ As above, $\left. \lambda
_{2}^{E_{11}}\right\vert _{T_{4}^{+}}\neq 0,$ which implies that $\lambda
_{2}^{E_{11}}$ does not change its sign when $\left( \mu _{1},\mu
_{2}\right) $ crosses $T_{4}^{+}.$ Thus, $E_{11}$ remains a repeller on $%
Y_{+}\cup R_{10}^{-}.$ Notice that 
\begin{equation*}
\lambda _{2}^{E_{11}}\left( \mu _{1},0\right) =\frac{1}{N\gamma }\sqrt{-\mu
_{1}N}\left(1+O\left(\mu_{1}\right)\right)>0\text{ if }\mu _{1}<0.
\end{equation*}

On the other hand, $\lambda _{2}^{E_{11}}$ changes its sign when $\left( \mu
_{1},\mu _{2}\right) $ crosses $T_{4}$ because $\left. \lambda
_{2}^{E_{11}}\right\vert _{T_{4}}=0.$ More exactly, $\lambda _{2}^{E_{11}}<0$
and $\lambda _{1}^{E_{11}}>0$ on $R_{10}^{+}.$ Thus, $E_{11}$ is a saddle on 
$R_{10}^{+}.$ $\square $

\bigskip

\begin{figure}[htpb]
\begin{center}
\includegraphics[width=0.3\textwidth]{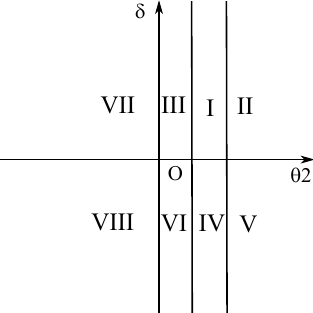}
\end{center}
\caption{When $\protect\delta \neq 0$ and $\protect\theta =0,$ eight cases in the $%
\protect\delta\protect\theta_{2}$- plane lead to eight bifurcation diagrams.}
\label{fig9}
\end{figure}

\begin{figure}[htpb]
\begin{center}
\includegraphics[width=1\textwidth]{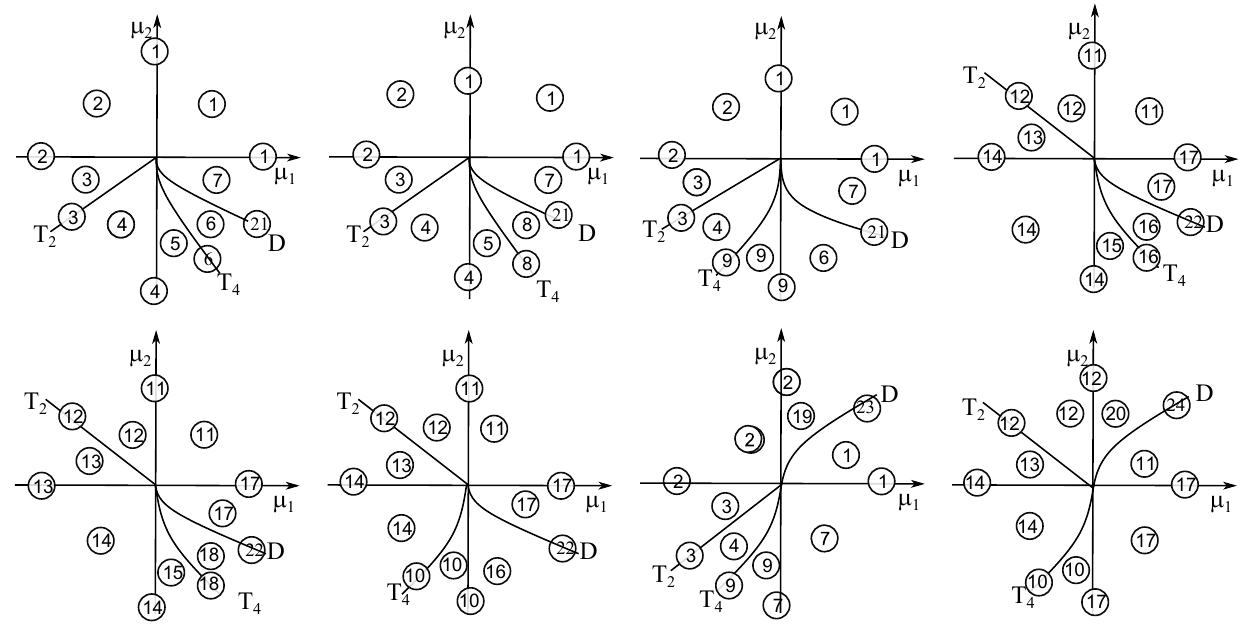}
\end{center}
\caption{Bifurcation diagrams for $\protect\delta \neq 0$ and $\protect%
\theta =0$ corresponding to the eight cases I-VIII.}
\label{fig10}
\end{figure}

\begin{figure}[htpb]
\begin{center}
\includegraphics[width=0.9\textwidth]{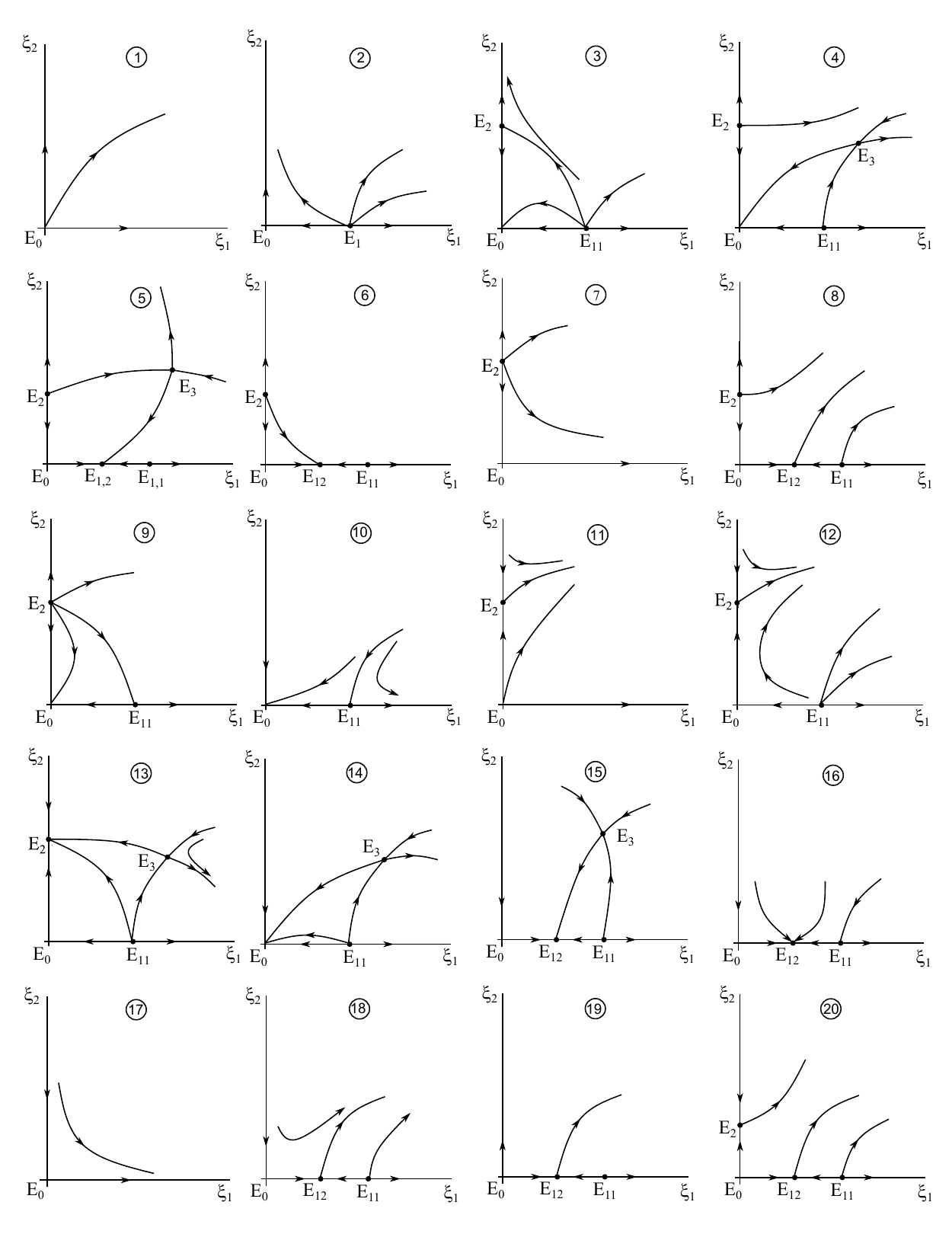}
\end{center}
\caption{Phase portraits corresponding to the bifurcation diagrams I-VIII,
when $\protect\delta \neq 0$ and $\protect\theta =0.$}
\label{fig11}
\end{figure}

Eight different cases arise in terms of $\delta $ and $\theta _{2},$ (Fig.%
\ref{fig9}), each one leading to a bifurcation diagram. There are 20
different regions in the eight bifurcation diagrams (Fig.\ref{fig10}). In
Table \ref{tabel4} we summarized the type of each equilibrium point from the
20 regions. The phase portraits corresponding to the 20 regions are depicted
in Fig.\ref{fig11}. The phase portraits on the saddle-node curve is
displayed in Fig.\ref{fig.12}.

\bigskip

I: $\delta >0,$ $\theta_{2}>0,$ $\theta _{2}-N\gamma>0$ and $%
\theta_{2}-2N\gamma<0,$ Fig.\ref{fig9} (I)

II: $\delta >0,$ $\theta_{2}>0,$ $\theta _{2}-N\gamma>0$ and $%
\theta_{2}-2N\gamma>0,$ Fig.\ref{fig9} (II)

III: $\delta >0,$ $\theta_{2}>0,$ $\theta_{2}-N\gamma<0$ and $%
\theta_{2}-2N\gamma<0,$ Fig.\ref{fig9} (III)

IV: $\delta <0,$ $\theta_{2}>0,$ $\theta_{2}-N\gamma>0$ and $%
\theta_{2}-2N\gamma<0,$ Fig.\ref{fig9} (IV)

V: $\delta <0,$ $\theta_{2}>0,$ $\theta_{2}-N\gamma>0$ and $%
\theta_{2}-2N\gamma>0,$ Fig.\ref{fig9} (V)

VI: $\delta <0,$ $\theta_{2}>0,$ $\theta_{2}-N\gamma<0$ and $%
\theta_{2}-2N\gamma<0,$ Fig.\ref{fig9} (VI)

VII: $\delta >0,$ $\theta_{2}<0,$ Fig.\ref{fig9} (VII)

VIII: $\delta <0,$ $\theta_{2}<0,$ Fig.\ref{fig9} (VIII)

\begin{table}[htpb]
\begin{tabular}{l|llllllllllllllllllll}
& $1$ & $2$ & $3$ & $4$ & $5$ & $6$ & $7$ & $8$ & $9$ & $10$ & $11$ & $12$ & 
$13$ & $14$ & $15$ & $16$ & $17$ & $18$ & $19$ & $20$ \\ \hline
\centering $O$ & $r$ & $s$ & $a$ & $a$ & $s$ & $s$ & $s$ & $s$ & $a$ & $a$ & 
$r$ & $s$ & $s$ & $a$ & $s$ & $s$ & $s$ & $s$ & $r$ & $r$ \\ 
$E_{11}$ & $-$ & $r$ & $r$ & $r$ & $r$ & $s$ & $-$ & $r$ & $s$ & $s$ & $-$ & 
$r$ & $r$ & $r$ & $r$ & $s$ & $-$ & $r$ & $r$ & $r$ \\ 
$E_{12}$ & $-$ & $-$ & $-$ & $-$ & $a$ & $a$ & $-$ & $s$ & $-$ & $-$ & $-$ & 
$-$ & $-$ & $-$ & $a$ & $a$ & $-$ & $s$ & $s$ & $s$ \\ 
$E_{2}$ & $-$ & $-$ & $s$ & $r$ & $r$ & $r$ & $r$ & $r$ & $r$ & $-$ & $s$ & $%
s$ & $a$ & $-$ & $-$ & $-$ & $-$ & $-$ & $-$ & $s$ \\ 
$E_{3}$ & $-$ & $-$ & $-$ & $s$ & $s$ & $-$ & $-$ & $-$ & $-$ & $-$ & $-$ & $%
-$ & $s$ & $s$ & $s$ & $-$ & $-$ & $-$ & $-$ & $-$ \\ 
&  &  &  &  &  &  &  &  &  &  &  &  &  &  &  &  &  &  &  & 
\end{tabular}%
\caption{\textit{The types of the equilibrium points of system \eqref{gl8}
for $\protect\delta \neq 0$ and $\protect\theta =0$ on different regions of
the bifurcation diagrams I-VIII; the abbreviations s, a, r stand for saddle,
attractor, repeller, respectively.}}
\label{tabel4}
\end{table}

\begin{figure}[htpb]
\begin{center}
\includegraphics[width=1\textwidth]{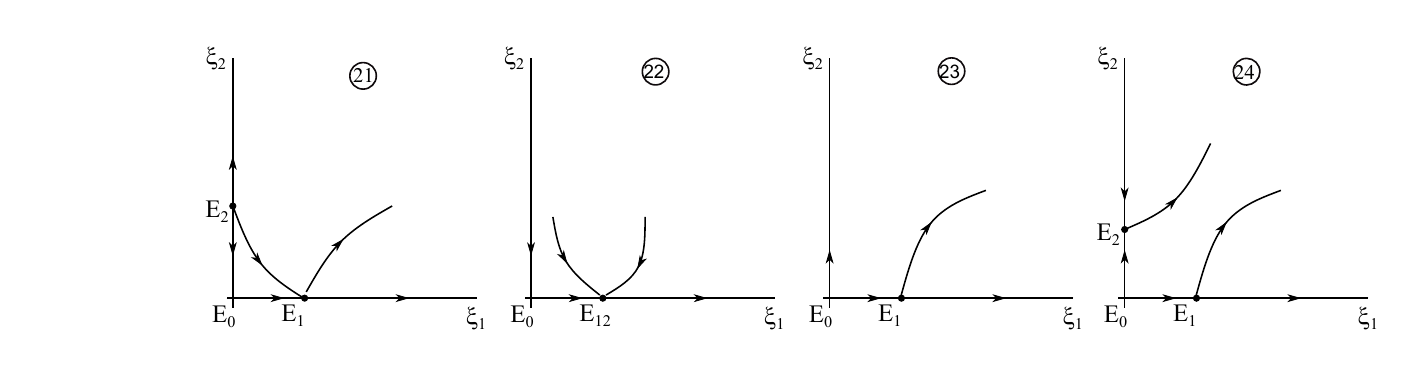}
\end{center}
\caption{Phase portraits corresponding to the saddle node curve, when $%
\protect\delta \neq 0$ and $\protect\theta =0.$}
\label{fig12}
\end{figure}

\begin{theorem}
Assume $\gamma _{1}\theta _{2}\delta \left( \theta _{2}-N\gamma \right) \neq
0,$ $\gamma >0$ and $\theta _{2}-2N\gamma \neq 0,$ where $\gamma _{1}=\frac{%
\partial \gamma \left( 0\right) }{\partial \mu _{1}}.$ Then $T_{4}$ is a
transcritical bifurcation curve.
\end{theorem}

\textit{Proof.} Assume first $\theta _{2}-2N\gamma <0,$ thus, $E_{3}$
collides to $E_{11}\left( -\gamma \mu _{2},0\right) $ on $T_{4}.$ Let $\mu
_{1}$ be the bifurcation parameter while $\mu _{2}<0$ is assumed fixed.
Denote by $\mu _{0}=\left( \mu _{1},\mu _{2}\right) \in T_{4},$ that is, $%
\mu _{1}=\gamma \left( \theta _{2}-N\gamma \right) \mu _{2}^{2}.$

On $T_{4},$ $E_{3}\left( \xi _{1},0\right) $ satisfies $\mu _{1}+\theta
\left( \mu _{0}\right) \xi _{1}+N\left( \mu _{0}\right) \xi _{1}^{2}=0$ and $%
\mu _{2}+\frac{1}{\gamma \left( \mu _{0}\right) }\xi
_{1}+R\left(\mu_{0}\right)\xi_{1}^{2}=0,$ thus, $\xi_{1}=\xi _{11}=-\gamma
\left( \mu _{0}\right) \mu _{2}.$ Whenever $\left( \mu _{1},\mu _{2}\right)
\in T_{4},$ denote by $\xi _{0}=\left( \xi _{11},0\right) .$

The Jacobian matrix at $\left( \xi _{0},\mu _{0}\right) $ is of the form 
\begin{equation*}
A=Df\left( \xi _{0},\mu _{0}\right) =\allowbreak \left( 
\begin{array}{cc}
\left( \theta \left( \mu _{0}\right) +2N\left( \mu _{0}\right) \xi
_{1}\right) \xi _{1} & \left( \gamma \left( \mu _{0}\right) +M\left( \mu
_{0}\right) \xi _{1}\right) \xi _{1} \\ 
0 & 0%
\end{array}%
\right) ,
\end{equation*}%
which has the eigenvalue $0$ with the eigenvector $v=\left( 
\begin{array}{cc}
v_{1} & 1%
\end{array}%
\right) ^{T},$ where $v_{1}=-\frac{\gamma \left( \mu _{0}\right) +M\left(
\mu _{0}\right) \xi _{1}}{\theta \left( \mu _{0}\right) +2N\left( \mu
_{0}\right) \xi _{1}}.$ The eigenvector corresponding to the eigenvalue $0$
in $A^{T}$ is $w=\allowbreak \left( 
\begin{array}{cc}
0 & 1%
\end{array}%
\right) ^{T}.$ Denote as above $f=\left( f_{1},f_{2}\right) .$

Further, $f_{\mu _{1}}$ has the form

\begin{equation*}
f_{\mu _{1}}=\left( 
\begin{array}{c}
\frac{\partial f_{1}}{\partial \mu _{1}} \\ 
\frac{\partial f_{2}}{\partial \mu _{1}}%
\end{array}%
\right) =\left( 
\begin{array}{c}
\xi _{1}\left( 1+\frac{\partial \theta \left( \mu \right) }{\partial \mu _{1}%
}\xi _{1}+\frac{\partial \gamma \left( \mu \right) }{\partial \mu _{1}}\xi
_{2}+\frac{\partial M\left( \mu \right) }{\partial \mu _{1}}\xi _{1}\xi _{2}+%
\frac{\partial N\left( \mu \right) }{\partial \mu _{1}}\xi _{1}^{2}+\frac{%
\partial L\left(\mu\right)}{\partial \mu_{1}}\xi_{2}^{2}\right) \\ 
\xi _{2}\left( G\left( \mu \right) \xi _{1}+\frac{\partial \delta \left( \mu
\right) }{\partial \mu _{1}}\xi _{2}+\frac{\partial S\left( \mu \right) }{%
\partial \mu _{1}}\xi _{1}\xi _{2}+\frac{\partial P\left( \mu \right) }{%
\partial \mu _{1}}\xi _{2}^{2}+\frac{\partial R\left(\mu\right)}{\partial
\mu_{1}}\xi_{1}^{2}\right)%
\end{array}%
\right) ,
\end{equation*}%
where $G\left( \mu \right) =-\frac{1}{\gamma ^{2}\left( \mu \right) }\frac{%
\partial \gamma \left( \mu \right) }{\partial \mu _{1}}.$ Then $%
C_{1}=w^{T}f_{\mu _{1}}\left( \xi _{0},\mu _{0}\right) =0.$ The Jacobian $%
Df_{\mu _{1}}$ at the point $\left( \xi _{0},\mu _{0}\right) $ applied to
the vector $v$ has the form

\begin{equation*}
Df_{\mu _{1}}\left( \xi _{0},\mu _{0}\right) =\left( 
\begin{array}{cc}
\frac{\partial ^{2}f_{1}}{\partial \xi _{1}\partial \mu _{1}}\left( \xi
_{0},\mu _{0}\right) & \frac{\partial ^{2}f_{1}}{\partial \xi _{2}\partial
\mu _{1}}\left( \xi _{0},\mu _{0}\right) \\ 
0 & G\left( \mu _{0}\right) \xi _{1}%
\end{array}%
\allowbreak \right) \left( 
\begin{array}{c}
v_{1} \\ 
1%
\end{array}%
\right) =\left( 
\begin{array}{c}
K \\ 
G\left( \mu _{0}\right) \xi _{1}%
\end{array}%
\right) ,
\end{equation*}%
where $K$ is an expression which is not needed in what follows. Then 
\begin{equation*}
C_{2}=w^{T}\left[ Df_{\mu _{1}}\left( \xi _{0},\mu _{0}\right) \left(
v\right) \right] =G\left( \mu _{0}\right) \xi _{1}=\frac{\gamma _{1}}{\gamma 
}\mu _{2}\left( 1+O\left( \mu _{2}\right) \right) \neq 0.
\end{equation*}%
Finally, we need to find $C_{3}=w^{T}\left[ D^{2}f\left( \xi _{0},\mu
_{0}\right) \left( v,v\right) \right] ,$ where $D^{2}f\left( \xi ,\mu
\right) \left( v,v\right) =\left( 
\begin{array}{c}
d^{2}f_{1}\left( \xi ,\mu \right) \left( v,v\right) \\ 
d^{2}f_{2}\left( \xi ,\mu \right) \left( v,v\right)%
\end{array}%
\right) .$ Since $w=\left( 
\begin{array}{cc}
0 & 1%
\end{array}%
\right) ^{T},$ we need to determine only $d^{2}f_{2}\left( \xi ,\mu \right)
\left( v,v\right) ,$ where $v=\left( v_{1},1\right) .$ These lead to 
\begin{equation*}
C_{3}=\allowbreak \allowbreak \frac{2}{\left( 2N\gamma -\theta _{2}\right)
\mu _{2}}\left( 1+O\left( \mu _{2}\right) \right) \allowbreak \neq 0,
\end{equation*}%
which confirms the claim. For $\theta _{2}-2N\gamma >0$ the proof is
similar. $\square $

\begin{remark}
One can show that $T_{2},$ $X^{+}=\left\{ \left( \mu _{1},0\right)
\left\vert \mu _{1}>0\right. \right\} ,$ $X^{-}=\left\{ \left( \mu
_{1},0\right) \left\vert \mu _{1}<0\right. \right\} ,$ $Y^{+}=\left\{ \left(
0,\mu _{2}\right) \left\vert \mu _{2}>0\right. \right\} $ and $Y^{-}=\left\{
\left( 0,\mu _{2}\right) \left\vert \mu _{2}<0\right. \right\} $ are also
transcritical bifurcation curves. The behavior of the system (\ref{gl8}) on the axes $X^{\pm},\ Y^{\pm}$ and on the transcritical curves $T_{2,4}$ is similar with the first degenerate case. On the saddle-node curves, the corresponding dynamics is presented
in Figure \ref{fig12}.
\end{remark}

\section{Conclusions}

In this paper we studied a generalized Lotka-Volterra model with small birth
rates of predator and pray. Three different cases have been considered: one
non-degenerate corresponding to $\delta \left( 0\right) \theta \left(
0\right) \neq 0,$ and two degenerate with $\theta (0)\neq 0$ and $\delta
(0)=0,$ respectively, $\theta (0)=0$ and $\delta (0)\neq 0.$

For the non-degenerate case, six different bifurcation diagrams emerged for
the description of the model's dynamics. The diagrams contain thirty
different regions in the parametric plane $\mu _{1}\mu _{2}.$ The type of
equilibria in each region and the corresponding phase portraits have been
obtained.

The model's dynamics in the two degenerate cases has been described by 16
bifurcation diagrams with 40 different regions. The phase portraits
corresponding to these regions have been presented. We showed that the
equilibrium points different from $O$ in the two degenerate cases, are born
or enter the first quadrant $Q_{1}$ mainly by saddle-node and transcritical
bifurcations. 

The model introduced in this work is far from being completely described. Many other
degeneracies may arise, such as $\theta \left( 0\right) =0$ and $\delta
\left( 0\right) =0,$ or $\theta_{1,2} =0$ and $\delta_{1,2} =0.$ These cases
remain open for further studies.

\section{Data statement} Data sharing not applicable to this article as no datasets were generated or analysed during the current study.

\section{Acknowledgments}

This research was supported by Horizon2020-2017-RISE-777911 project.

\end{document}